\documentclass[12pt]{amsart}
\usepackage{amsmath,amssymb,graphicx}
\usepackage[all]{xy}
\usepackage[bbgreekl]{mathbbol}
\usepackage{verbatim}

\newtheorem{thm}{Theorem}[subsection]

\newtheorem{lem}[thm]{Lemma}
\newtheorem{prop}[thm]{Proposition}
\theoremstyle{definition}
\newtheorem{defn}[thm]{Definition}
\theoremstyle{remark}
\newtheorem{rem}[thm]{Remark}
\numberwithin{equation}{section}

\DeclareFontFamily{U}{rsf}{} \DeclareFontShape{U}{rsf}{m}{n}{
  <5> <6> rsfs5 <7> <8> <9> rsfs7 <10->  rsfs10}{}
\DeclareMathAlphabet{\mathscr}{U}{rsf}{m}{n}

\newcommand{\mycal}[1]{\mathscr{#1}}

\newcommand{\mto}[1]{\stackrel{#1}{\longrightarrow}}
\newcommand{\abs}[1]{\left\vert#1\right\vert}

\newcommand{\cx}{\mathbb C}

\newcommand{\A}{\mathcal{A}}

\newcommand{\Cc}{\,{\mathcal C}}
\newcommand{\Pc}{{\mathcal P}}
\newcommand{\Pcc}{{\mycal P}}

\newcommand{\Zb}{{\mathbb Z}}

\newcommand{\V}{\mathcal{V}}
\newcommand{\D}{\mathcal{D}}

\newcommand{\As}{\text{\bfseries\sf{A}}}

\newcommand{\Ec}{\mathcal E}

\newcommand{\Al}{\mathfrak a}

\renewcommand{\imath}{\sqrt{-1}}
\newcommand{\Sc}{\mycal S}

\newcommand{\Hc}{\mathcal H}

\newcommand{\Eb}{\mathbb E}
\newcommand{\Fb}{\mathbb F}

\newcommand{\cf}{2\pi\imath\,}

\newcommand{\Dc}{\mathfrak{D}}

\newcommand{\np}{ }
\newcommand{\dbar}{\overline{\partial}}
\DeclareMathOperator{\modl}{Mod} 
\DeclareMathOperator{\op}{\circ} 
\DeclareMathOperator{\Hom}{Hom}
\DeclareMathOperator{\Ext}{Ext} 
\DeclareMathOperator{\End}{End}
\DeclareMathOperator{\Ob}{Ob} 
\DeclareMathOperator{\Ho}{Ho}

\DeclareMathOperator{\I}{Im}

\begin{document}

\title[Duality and equivalence]{Duality and equivalence of module categories in noncommutative geometry}
\author{Jonathan Block}\thanks{J.B. partially supported by NSF grant DMS02-04558}%
\address{Department of Mathematics, University of Pennsylvania, Philadelphia, PA 19104}%
\email{blockj@math.upenn.edu}%
\dedicatory{In memory of Raoul}
\thanks{}

\begin{abstract}

\noindent  We develop a general framework 
to describe dualities from algebraic, differential, and noncommutative geometry, as well as physics. We pursue a relationship between the
Baum-Connes  conjecture in operator $K$-theory  and  derived
equivalence statements in algebraic geometry and physics. We associate to certain data, reminiscent of spectral triple data, a
differential graded category in such a way that we can recover the
derived category of coherent sheaves on a complex manifold. 
\end{abstract}
\maketitle 
\section*{Introduction}  
In various geometric contexts, there are
duality statements that are expressed in terms of appropriate
categories of modules. We have in mind, for example, 
the Baum-Connes conjecture from noncommutative geometry, T-duality and Mirror symmetry from 
complex geometry and mathematical physics. This is the first in a series of papers that sets up a framework to study and unify these dualities from a noncommutative geometric point of view. We also view this project as an attempt to connect the noncommutative geometry of Connes, \cite{Co1} with the categorical approach to noncommutative geometry, represented for example by Manin and Kontsevich. 

Traditionally, the complex structure is encoded in the sheaf of
holomorphic functions. However, for situations we have in mind
coming from noncommutative geometry, one can not use local types
of constructions, and we are left only with global differential
geometric ones. A convenient setting to talk about integrability of geometric
structures and the integrability of geometric structures on their
modules is that of a differential graded algebra and more generally, a curved differential graded algebra. Thus, for
example, a complex structure on a manifold is encoded in its
Dolbeault algebra $\As=(\A^{0,\bullet}(X),\dbar)$, and a holomorphic vector bundle can be viewed as the data of a 
finitely generated projective module over $\A^{0,0}$ together with a flat $\dbar$-connection. Similarly, holomorphic gerbes can be encoded in terms of a curved differential graded algebra with non-trivial curvature. Curved dgas  appear naturally in the context of matrix factorizations and Laundau-Ginzberg models, \cite{Ei}, \cite{Or}. Indeed, these fit very easily into our framework. 

Of course, one is interested in more modules than just the
finitely generated projective ones. In algebraic geometry, the
notion of coherent module is fundamental. In contrast to
projective algebraic geometry however, not every coherent sheaf
has a resolution by vector bundles; they only locally have such
resolutions. Toledo and Tong, \cite{TT1}, \cite{TT2}, handled this
issue by introducing twisted complexes. Our construction is a
global differential geometric version of theirs.

We have found the language of
differential graded categories to be useful, \cite{BK}, \cite{Ke}, \cite{Dr}.
In particular, for a curved dga $\As$ we construct a very natural differential
graded category $\Pc_\As$ which can then be derived. The desiderata of such
a category are
\newline
\begin{itemize}
    \item it should be large enough
to contain in a natural way the coherent holomorphic sheaves (in
the case of the Dolbeault algebra), and
    \item it should be flexible enough to
allow for some of Grothendieck's six operations, so that we can
prove Mukai type duality statements.
\end{itemize}

The reason for introducing $\Pc_\As$ is that the ordinary category of dg-modules over the Dolbeault dga has the wrong homological algebra; it has the wrong notion of quasi-isomorphism. A morphism between complexes of holomorphic vector bundles considered as dg-modules over the Dolbeault algebra is a quasi-isomorphism if it induces an isomorphism on the total complex formed by the gloabal sections of  the Dolbeault algebra with values in the complexes of holomorphic vector bundles, which is isomorphic to  their hypercohomology. On the other hand, $\Pc_\As$ and the modules over it have the correct notion of quasi-isomorphism. In particular, $\Pc_\As$ is not an invariant of quasi-isomorphism of dga's. To be sure, we
would not want this. For example, the dga which is $\cx$ in degree
$0$ and $0$ otherwise is quasi-isomorphic to the Dolbeault algebra
of $\mathbb{CP}^n$. But $\mathbb{CP}^n$ has a much richer module
category than anything $\cx$ could provide. We show that the homotopy category of $\Pc_\As$ where $\As$ is the Dolbeault algebra of a compact complex manifold $X$ is equivalent to the derived category of sheaves of $\mathcal{O}_X$-modules with coherent cohomology. Our description of the coherent derived category has recently been used by Bergman, \cite{Be} as models for $B$-branes.

To some extent, what we do is a synthesis of Kasparov's
$KK$-theory, \cite{Ka} and of Toledo and Tong's twisted complexes,
\cite{TT1}, \cite{TT2}, \cite{OTT}.

\subsection*{In appreciation of Raoul Bott} I am always amazed by the profound impact that he had, and still has on my life. During the time I was his student, I learned much more from him than mere mathematics. It was his huge personality, his magnanimous heart, his joy in life and his keen aesthetic that has had such a lasting affect. I miss him.  
\subsection*{Acknowledgements.} We would like to thank Oren Ben-Bassat, Andre Caldararu, Calder Daenzer, Nigel
Higson, Anton Kapustin, the referee, Steve Shnider, Betrand Toen and especially Tony Pantev
for many conversations and much guidance regarding this project.

\section{Baum-Connes and Fourier-Mukai} There are two major motivations for our project. The first is to have a general framework that will be useful in dealing with categories of modules that arise in geometry and physics. For example, we will apply our framework to construct categories of modules over symplectic manifolds. Second, as mentioned earlier in the introduction, this series of papers is
meant to pursue a relationship between \begin{enumerate} \item the
Baum-Connes conjecture in operator $K$-theory  and \item derived
equivalence statements in algebraic geometry and physics.
\end{enumerate}
In particular, we plan to refine, in certain cases, the
Baum-Connes conjecture from a statement about isomorphism of two
topological $K$-groups to a derived equivalence of categories
consisting of modules with geometric structures, for example, coherent sheaves on complex manifolds.  We will see that
there are natural noncommutative geometric spaces that are derived equivalent to classical
algebraic geometric objects.

Let us explain the obvious formal analogies between (1) and (2).
For simplicity let $\Gamma$ be a discrete torsion free group with
compact $B\Gamma$. In this situation, the Baum-Connes conjecture
says that an explicit map, called the assembly map,
\begin{equation}\label{BC} \mu:K_*(B\Gamma)\to K_*(C^*_r\Gamma)
\end{equation}
is an isomorphism. Here $C^*_r\Gamma$ denotes the reduced group
$C^*$-algebra of $\Gamma$. The assembly map can be described in
the following way. On $C(B\Gamma)\otimes C^*_r\Gamma $ there is a
finitely generated projective right module $\Pcc$ which can be
defined as the sections of the bundle of $C^*_r\Gamma$-modules
$$E\Gamma\times_\Gamma C^*_r\Gamma.$$ 
This projective module is a ``line bundle" over $C(B\Gamma)\otimes
C^*_r\Gamma$. Here, $C(X)$ denote the complex valued continuous functions on a compact space $X$. The assembly map is the map defined by taking the
Kasparov product with $\Pcc$ over $C(B\Gamma)$. This is some sort
of index map.
$$\mu:x\in KK(C(B\Gamma), \mathbb{C})\longmapsto x\cup \Pcc \in KK(\mathbb{C},C^*_r\Gamma)$$ where $\Pcc \in KK(\mathbb{C},C(B\Gamma)\otimes
C^*_r\Gamma)$

We now describe Mukai duality in a way that makes it clear that it
refines Baum-Connes.
Now let $X$ be a complex torus. Thus $X=V/\Lambda$ where $V$ is a
$g$-dimensional complex vector space and  $\Lambda\cong \Zb^{2g}$
is a lattice in $V$. Let $X^\vee$ denote the dual complex torus.
This can be described in a number of ways:

\begin{itemize}
\item as $\mbox{Pic}^0(X)$, the manifold of holomorphic line
bundles on $X$ with first Chern class $0$ (i.e., they are
topologically trivial); \item as the moduli space of flat unitary
line bundles on $X$. This  is the same as the space of irreducible unitary
representations of $\pi_1(X)$, but it has a complex structure that
depends on that of $X$; \item and most explicitly as
$\overline{V}^\vee/\Lambda^\vee$ where $\Lambda^\vee$ is the dual
lattice,
$$\Lambda^\vee=\{v\in\overline{V}^\vee \, |\, \mbox{ Im } <v,\lambda>\in
\Zb\,\, \forall \lambda \in\Lambda\}.$$ Here $\overline{V}^\vee$
consists of conjugate linear homomorphisms from $V$ to
$\mathbb{C}$.
\end{itemize}

We note that $X=B\Lambda$ and that $C(X^\vee)$ is canonically
$C^*_r\Lambda$. Hence Baum-Connes predicts (and in fact it is
classical in this case) that $K_*(X)\cong K_*(C^*_r\Lambda)\cong K^*(X^\vee)$.

On $X\times X^\vee$ there is a canonical line bundle, $\Pcc$, the
Poincar\'{e} bundle which is uniquely determined by the following
universal  properties:
\begin{itemize}
\item $\Pcc|X\times\{p\}\cong p$ where $p\in X^\vee$ and is
therefore a line bundle on $X$. \item $\Pcc|\{0\}\times X^\vee$ is
trivial.
\end{itemize}
Now Mukai duality says that there is an equivalence of derived
categories of coherent sheaves
$$D^b(X)\to D^b(X^\vee)$$
induced by the functor
$${\mathcal F}\mapsto p_{2*}(p_1^*{\mathcal F}\otimes \Pcc)$$
where $p_i$ are the two obvious projections. The induced map at
the level of $K_0$ is an isomorphism and is clearly a holomorphic version of 
the Baum-Connes
Conjecture for the group $\Lambda$.

\section{The dg-category $\Pc_\As$ of a curved dga}
\subsection{dg-categories}

\begin{defn} For complete definitions and facts regarding dg-categories, see \cite{BK}, \cite{Dr}, \cite{Ke} and \cite{Ke2}. Fix a field $k$. A {\em differential graded category (dg-category) } is a
category enriched over $\Zb$-graded complexes (over $k$) with differentials increasing degree. That is,
a category $\Cc$ is a dg-category if for $x$ and $y$ in $\Ob\Cc$
the hom set
\[\Cc(x,y)\]
forms a $\Zb$-graded complex of $k$-vector spaces. Write
$(\Cc^\bullet(x,y),d)$ for this complex, if we need to reference the degree or differential in the complex. In addition, the
composition, for $x,y,z\in \Ob\Cc$
\[
\Cc(y,z)\otimes \Cc(x,y)\to
\Cc(x,z)\] is a morphism of complexes. Furthermore, there are obvious associativity and unit axioms. 
\end{defn}

\subsection{Curved dgas}
In many situations the integrability conditions are not expressed
in terms of flatness but are defined in terms of other curvature
conditions. This leads us to set up everything in the more general
setting of curved dga's. These are dga's where $d^2$ is not necessarily zero. 

\begin{defn} A curved dga \cite{PP} (Schwarz \cite{S} calls them $Q$-algebras) is a triple
\[\As=(\A^\bullet, d, c)\] where $\A^\bullet$ is a (non-negatively)
graded algebra over a field  $k$ of characteristic $0$, with a derivation
\[d:\A^\bullet\to \A^{\bullet+1}\]
which satisfies the usual graded Leibniz relation but
\[d^2(a)=[c,a]\]
where $c\in \A^2$ is a fixed element (the curvature). Furthermore
we require the Bianchi identity $dc=0$. Let us write $\A$ for the degree $0$ part of $\A^\bullet$, the ``functions" of $\As$. \end{defn} A dga is the
special case where $c=0$. Note that $c$ is part of the data and
even if $d^2=0$, that $c$ might not be $0$, and gives a non dga
example of a curved dga. The prototypical example of a curved dga
is $(\A^\bullet(M,\End(E)), Ad \nabla, F)$ of differential forms
on a manifold with values in the endomorphisms of a vector bundle
$E$ with connection ${\nabla}$ and curvature $F$.

\subsection{The dg-category $\Pc_{\As}$}
Our category $\Pc_\As$ consists of special types of $\As$-modules. We start with a $\Zb$-graded  right module $E^\bullet$ over $\A$.

\begin{defn}\label{Z-connection} A
$\Zb$-connection $\Eb$  is a
$k$-linear  map
\[\Eb: E^\bullet\otimes_\A\A^\bullet\to E^\bullet\otimes_\A\A^\bullet
\]
of total degree one, which satisfies the usual Leibniz condition
\[  \Eb(e\omega)=(\Eb(e\otimes 1))\omega+(-1)^e e d\omega\]
\end{defn}
Such a connection is determined by its value on $E^\bullet$. Let
$\Eb^k$ be the component of $\Eb$ such that $\Eb^k:E^\bullet\to
E^{\bullet-k+1}\otimes_\A\A^k$, thus
$\Eb=\Eb^0+\Eb^1+\Eb^2+\cdots$.  It is clear that $\Eb^1$ is a
connection on each component $E^n$ in the ordinary sense (or the negative of a connection if $n$ is odd) and that $\Eb^k$
is $\A$-linear for $k\ne 1$. 

\np Note that for a $\mathbb{Z}$-connection $\mathbb{E}$ on $E^\bullet$ over a curved dga $\As=(\A^\bullet, d, c)$,  the usual curvature $\mathbb{E}\circ\mathbb{E}$ is not $\A$-linear. Rather, we define the {\em relative curvature} to be the operator  
\[
F_\Eb(e)=\mathbb{E}\circ\mathbb{E}(e)+e\cdot c
\]
and this is $\A$-linear.

\begin{defn} For a curved dga $\As=(\A^\bullet,d,c)$, we define the dg-category $\Pc_{\As}$:
\begin{enumerate}
\item An object $E=(E^\bullet, \mathbb{E})$ in
$\Pc_{\As}$, which we call a {\em cohesive module}, is a
$\Zb$-graded (but bounded in both directions) right module
$E^\bullet$ over $\A$ which is finitely generated and projective,
together with a $\Zb$-connection
\[\Eb: E^\bullet\otimes_\A\A^\bullet\to E^\bullet\otimes_\A\A^\bullet
\]
that satisfies the integrability condition that the relative curvature vanishes
\[F_\mathbb{E}(e)=\Eb\circ\Eb(e)+e\cdot
c=0\]
for all $e\in E^\bullet$.

\item The morphisms of degree $k$,  ${\Pc^k_{\As}}(E_1,E_2)$ between two cohesive modules $E_1=(E_1^\bullet,
\Eb_1)$ and $E_2=(E_2^\bullet,\Eb_2)$ of degree $k$ are 
\[
\{\phi:E_1^\bullet\otimes_\A
\A^\bullet\to E_2^{\bullet}\otimes_\A \A^\bullet\,\,|\,\,\mbox{ of degree $k$ and }
\phi(ea)=\phi(e)a \,\,\, \forall a\in \A^\bullet\}
\]with differential defined in the standard way
\[
d(\phi)(e)= \Eb_2(\phi(e))-(-1)^{\abs{\phi}}\phi(\Eb_1(e))
\]
Again, such a $\phi$ is determined by its restriction to
$E_1^\bullet$ and if necessary we denote the component of $\phi$
that maps \begin{equation}\label{pcomp} E_1^\bullet\to
E_2^{\bullet+k-j}\otimes_\A \A^j \end{equation} by $\phi^j$. 
\end{enumerate}
\end{defn}
Thus ${\Pc^k_{\As}}(E_1,E_2)=\Hom^k_\A(E_1^\bullet,E_2^{\bullet}\otimes_\A \A^\bullet)$
\begin{prop} For  $\As=(\A^\bullet, d, c)$ a curved dga, the category $\Pc_{\As}$ is a dg-category.
\end{prop}
\np This is clear from the following lemma.

\begin{lem}
Let $E_1$, $E_2$ be cohesive modules over the curved dga
$\As=(\A^\bullet,d,c)$. Then the differential defined above
\[d:{\Pc^\bullet_{\As}}(E_1,E_2)\to {\Pc^{\bullet+1}_{\As}}(E_1,E_2)\]
satisfies $d^2=0$.
\end{lem}

\subsection{The homotopy category and triangulated structure}
Given a dg-category $\mathcal{C}$, one can form the subcategory $Z^0\mathcal{C}$ which has the same objects as $\mathcal{C}$ and whose morphisms  from an object $x\in \mathcal{C}$ to an object $y\in \mathcal{C}$ are the degree $0$ closed morphisms in $\Cc(x,y)$.  We also form the homotopy category $\Ho \mathcal{C}$ which has the same objects as $\mathcal{C}$ and whose morphisms  are the $0$th cohomology  
\[\Ho\mathcal{C}(x,y)=H^0(\mathcal{C}(x,y)).\]
 
We define a shift functor on the category $\Pc_{\As}$. For $E=(E^\bullet, \mathbb{E})$ set 
$E[1]=(E[1]^\bullet,\mathbb{E}[1])$ where $E[1]^\bullet=E^{\bullet+1}$ and $\mathbb{E}[1]=-\mathbb{E}$.
It is easy to verify that $E[1]\in \Pc_{\As}$. Next for $E, F\in \Pc_{\As}$ and $\phi\in Z^0\Pc_{\As}(E,F)$, define 
the cone of $\phi$, 
$\mbox{Cone}(\phi)=(\mbox{Cone}(\phi)^\bullet, \mathbb{C}_\phi)$ by
\[
\mbox{Cone}(\phi)^\bullet=\left(\begin{array}{c}F^\bullet \\  \oplus\\ E[1]^\bullet\end{array}\right)
\]
and 
\[
\mathbb{C}_\phi=\left(\begin{array}{cc} \mathbb{F} 		& \phi\\
                                 0              & \mathbb{E}[1] \end{array}\right)
\]
We then have a triangle of degree $0$ closed morphisms 
\begin{equation}\label{distinguishedtriangles}
E \stackrel{\phi}{\longrightarrow} F \longrightarrow \mbox{Cone}(\phi)    \longrightarrow E[1]  
\end{equation}

\begin{prop} Let $\As$ be a curved dga. Then the dg-category $\Pc_{\As}$ is pretriangulated in the sense of Bondal and Kapranov, \cite{BK}. Therefore, the category $\Ho\Pc_{\As}$ is triangulated with the collection of distinguished triangles being isomorphic to those of the form \ref{distinguishedtriangles}.
\end{prop} \begin{proof} The proof of this is the same as that of Proposition 1 and Proposition 2 of \cite{BK}.\end{proof}

\subsection{Homotopy equivalences}

\np As described above, a degree $0$ closed morphism $\phi$ between cohesive modules $E_i=(E_i^\bullet,\Eb_i)$, $i=1, 2$,
over $\As$ is a {\em homotopy equivalence} if it induces an isomorphism in $\Ho \Pc_{\As}$. We want to give a simple criterion for $\phi$ to define such a homotopy equivalence. On the complex ${\Pc_{\As}}(E_1,E_2)$
define a decreasing filtration  by
\[ F^k{\Pc^j_{\As}}(E_1,E_2)=\{\phi\in
{\Pc^j_{\As}}(E_1,E_2)| \,\,\,\phi^i=0 \mbox{ for }
i<k\}\] where $\phi^i$ is defined as in \eqref{pcomp}.

\begin{prop} There is a spectral sequence
\begin{equation}\label{SS} E_0^{pq}\Rightarrow
H^{p+q}({\Pc^\bullet_{\As}}(E_1,E_2))\end{equation}
where
\[ E_0^{pq}= \mbox{gr}{\Pc^\bullet_{\As}}
(E_1,E_2)=\{\phi^p\in{\Pc^{p+q}_{\As}}(E_1,E_2):E_1^\bullet\to
E_2^{\bullet+q}\otimes_\A\A^p\}\] with differential
\[d_0(\phi^p)=\Eb^0_2\circ\phi^p-(-1)^{p+q}\phi^p\circ\Eb_1^0\]
\end{prop}

\begin{prop}
A closed morphism $\phi\in
{\Pc^0_{\As}}(E_1,E_2)$ is a homotopy equivalence  if and only if
$\phi^0:(E_1^\bullet,\Eb_1^0)\to (E_2^\bullet,\Eb_2^0)$ is a
quasi-isomorphism of complexes of $\A$-modules.
\end{prop}
\begin{proof} Let $E=(E^\bullet, \Eb)$ be any object in
$\Pc_{\As}$. Then $\phi$ induces a map of complexes

\begin{equation}\label{a} \phi:\Pc_{\As}^\bullet(E,E_1)\to {\Pc^\bullet_{\As}}(E,E_2)\end{equation}
We show that  the
induced map on $E_1$-terms of the spectral sequences are
isomorphisms. Indeed, the quasi-isomorphism of $(E_1^\bullet,\Eb_1^0)\to (E_2^\bullet,\Eb_2^0)$ imply that they are actually chain homotopy equivalent since $E_1^\bullet$ and $E_2^\bullet$ are projective. Hence for each $p$ 
that
\[\phi^0\otimes {\mathbb I}:( E_1^\bullet\otimes_\A\A^p,\Eb^0_1\otimes\mathbb{I})
\to (E_2^\bullet\otimes_\A\A^p,\Eb^0_2\otimes\mathbb{I})\] is a
quasi-isomorphism and then
\[\mbox{gr }(\phi)=\phi^0:E_0^{pq}\cong\Hom_\A(E^\bullet,E_1^{\bullet+q}\otimes_\A\A^p)\to
E_0^{pq}\cong\Hom_\A(E^\bullet,E_2^{\bullet+q}\otimes_\A\A^p)\] is
a quasi-isomorphism after one last double complex argument since
the modules $E^\bullet$ are projective over $\A$ .  Thus \ref{a} is a quasi-isomorphism for all $E$ and this implies $\phi$ is an isomorphism in $\Ho\Pc_\As$. 

The other direction follows easily. \end{proof}

\subsection{The dual of a cohesive module}\label{dualtwistedcomplex}
 We define a duality functor
which will be of use in future sections.  Let $\As=(\A^\bullet, d, c)$
be a curved dga. Its opposite is $\As^{\op}=(\A^{\op\bullet}, d, -c)$
where $\A^{\op\bullet}$ is the graded algebra whose product is
given by
\[
a\cdot_{\op}b=(-1)^{|a||b|}ba
\]
We will not use the notation $\cdot_{\op}$ for the product any
longer. We can now define the category of left cohesive modules
over $\As$ as $\Pc_{\As^{\op}}$.

We define the duality dg functor
\[
^\vee:\Pc_{\As}\to \Pc_{\As^{\op}}
\]
by \[E=(E^\bullet,\Eb)\mapsto E^\vee=(E^{\vee\bullet},\Eb^\vee)
\]
where $E^{\vee k}=\Hom_\A(E^{-k},\A)$ and for $\phi\in
E^{\vee\bullet}$
\[
(\Eb^\vee\phi) (e)=d(\phi(e))-(-1)^{|\phi|}\phi(\Eb (e))
\]
It is straightforward that $E^\vee$ is indeed cohesive over $\As^{\op}$. There is a natural pairing between $E$ and $E^\vee$. And
moreover the connection was defined so that the relation
\[
\langle \Eb^\vee(\phi),e\rangle + (-1)^{|\phi|}\langle \phi,\Eb (e)\rangle =
d\langle \phi,e\rangle
\]
holds.
Note that the complex of morphisms $\Pc^\bullet_{\As}(E_1,E_2)$ between cohesive modules can be identified with 
\[
(E_2^\bullet\otimes_\A\A^\bullet\otimes_\A E_1^{\vee\bullet}, 1\otimes 1\otimes \Eb_1^\vee+1\otimes d\otimes 1+\Eb_2\otimes 1\otimes 1)
\]
\subsection{Functoriality}
We now discuss a construction of functors between categories of the
form $\Pc_{\As}$. Given two curved dga's,
$\As_1=(\A_1^\bullet,d_1,c_1)$ and $\As_2=(\A_2^\bullet,d_2,c_2)$ a
homomorphism from $\As_1$ to $\As_2$ is a pair
$(f,\omega)$ where $f:\A_1^\bullet\to\A_2^\bullet$ is a morphism
of graded algebras, $\omega\in \A_2^1$ and they satisfy
\begin{enumerate}
\item $f(d_1a_1)=d_2f(a_1)+[\omega,f(a_1)]$ and \item
$f(c_1)=c_2+d_2\omega+\omega^2$.
\end{enumerate}
Given a homomorphism of curved dga's $(f,\omega)$ we define a dg
functor
\[f_*:\Pc_{\As_1}\to \Pc_{\As_2}\]
as follows. Given $E=(E^\bullet,\Eb)$ a cohesive module over
$\As_1$, set $f^*(E)$ to be the cohesive module over
$\As_2$
\[(E^\bullet\otimes_{\A_1} \A_2, \Eb_2)\]
where $\Eb_2(e\otimes b)=\Eb(e)b+e\otimes(d_2b+\omega b)$. One
checks that $\Eb_2$ is still a $\Eb$-connection and satisfies
\[ (\Eb_2)^2(e\otimes b)=-(e\otimes b)c_2\]

\np This is a special case of the following construction. Consider
the following data, $X=(X^\bullet,\mathbb{X})$ where
\begin{enumerate}
\item $X^\bullet$ is a graded finitely generated projective
right-$\A_2$-module, \item $\mathbb{X}:X^\bullet\to
X^\bullet\otimes_{\A_2}\A_2^\bullet$ is a $\Zb$-connection, \item
$\A^\bullet_1$ acts on the left of
$X^\bullet\otimes_{\A_2}\A_2^\bullet$ satisfying

\[
a\cdot (x\cdot b)=(a\cdot x)\cdot b
\]
and
\[
\mathbb{X}(a\cdot(x\otimes b))=da\cdot(x\otimes
b)+a\cdot\mathbb{X}(x\otimes b)
\]
for $a\in \A_1^\bullet$, $x\in X^\bullet$ and $b\in \A_2^\bullet$,
 \item
$\mathbb{X}$ satisfies the following condition
\[
\mathbb{X}\circ\mathbb{X}(x\otimes b)=c_1\cdot (x\otimes
b)-(x\otimes b)\cdot c_2
\]
 on the complex
$X^\bullet\otimes_{\A_2}\A_2^\bullet$.

\end{enumerate}
Let us call such a pair $X=(X^\bullet,\mathbb{X})$ an
$\As_1-\As_2$-cohesive bimodule.

\np Given an $\As_1-\As_2$-cohesive bimodule
$X=(X^\bullet,\mathbb{X})$, we can then define a dg-functor (see the next section for the definition)
\[ X^*:\Pc_{\As_1}\to \Pc_{\As_2}\]
by 
\[X^*(E^\bullet, \Eb)=(E^\bullet\otimes_{\A_1}X^\bullet,\Eb_2)\]
where $\Eb_2(e\otimes x)=\Eb(e)\cdot x+e\otimes \mathbb{X}(x)$,
where the $\cdot$ denotes the action of $\A_1^\bullet$ on
$X^\bullet\otimes_{\A_1}\A_2^\bullet$.  One easily checks that
$X^*(E)$ is an object of $\Pc_{\As_2}$. We will write
$\mathbb{E}\otimes \mathbb{X}$ for $\mathbb{E}_2$.

\begin{rem}
\begin{enumerate} \item The previous case of a homomorphism between
curved dga's occurs by setting $X^\bullet=\A_2$ in degree $0$.
$\A_1^\bullet $ acts by the morphism $f$ and the $\Zb$-connection
is
\[ \mathbb{X}(a_2)=d_2(a_2)+\omega\cdot a_2.\]

\item To give another example of an
$\As_1-\As_2$-cohesive bimodule, consider a manifold
$M$ with two vector bundles with connection $(E_1,\nabla_2)$ and
$(E_2,\nabla_2)$. Let $c_i$ be the curvature of $\nabla_i$. Set
$\As_i=(\A_i^\bullet, d_i,c_i)=(\A^\bullet(M; \End(E_i), \mbox{Ad}\nabla_i)$. 
Then we define a cohesive  bimodule between
them by setting
\[
X^\bullet=\Gamma(M;\Hom(E_2,E_1))
\]
in degree $0$. $X^\bullet$ has a $\Zb$-connection

\[\mathbb{X}(\phi)(e_2)=\nabla_1(\phi(e_2))-\phi(\nabla_2
e_2)\] and maps $X^\bullet\to
X^\bullet\otimes_{\A_2}\A_2^\bullet$. Then
$(\mathbb{X})^2(\phi)=c_1\cdot \phi-\phi\cdot c_2$ as is required.
This cohesive bimodule implements a dg-quasi-equivalence between
$\Pc_{\As_1}$ and $\Pc_{\As_2}$. (See the next section for the definition of a dg-quasi-equivalence.)

\end{enumerate}
\end{rem}

\section{Modules over $\Pc_\As$}
It will be important for us to work with modules over $\Pc_{\As}$ and not just with the objects of $\Pc_{\As}$ itself. 
\subsection{Modules over a dg-category} We first collect some general definitions, see \cite{Ke} for more details. 
\begin{defn} A functor $F:\Cc_1\to \Cc_2$ between two dg-categories is a {\em
dg-functor} if the map on hom sets
\begin{equation}\label{dg-functor}F:{\Cc_1}(x,y)\to
{\Cc_2}(Fx,Fy)\end{equation}  is a chain map of
complexes. A dg-functor $F$ as above is a {\em quasi-equivalence}
if the maps in (\ref{dg-functor}) are quasi-isomorphisms and
$\Ho(F):\Ho \Cc_1\to \Ho\Cc_2$ is an equivalence of categories.
\end{defn}

\np Given a dg-category $\Cc$ , one can define the category of
(right) dg-modules over $\Cc$, $\modl-\Cc$. This
consists of dg-functors from the opposite dg-category $\Cc^\circ$
to the dg-category  $\Cc(k)$ of complexes over $k$. More explicitly, a right $\Cc$-module $M$ is an assignment to each 
$x\in\Cc$, a complex $M(x)$ and chain maps for any $x, y\in\Cc$ 
\begin{equation}\label{action}
M(x)\otimes\Cc(y,x)\to M(y)
\end{equation}
satisfying the obvious associativity and unit conditions. A morphism $f\in\modl-\Cc(M,N)$ between right $\Cc$-modules $M$ and $N$  is 
an assignment of a map of complexes 
\begin{equation}\label{modmorph}
f_x: M(x)\to N(x)
\end{equation}
for each object $x\in \Cc$ compatible with the maps in \eqref{action}. Such a map is called a {\em quasi-isomorphism}
if $f_x$ in \eqref{modmorph} is a quasi-isomorphism of complexes for each $x\in\Cc$. 
One can make modules over a dg-category into a dg-category itself. The morphisms we have defined in $\modl-\Cc$ are the degree $0$ closed morphisms of this dg-category. The category of left modules $\Cc-\modl$ is defined in an analogous way. 

The category  $\modl-\Cc$ has a model structure used by Keller to define its derived category, \cite{Ke}, \cite{Ke2}. The quasi-isomorphisms in $\mbox{Mod-}\Cc$ are those we just defined. The fibrations are the component-wise surjections and the cofibrations are defined by the usual lifting property. Using this model structure we may form the homotopy category of $\mbox{Mod-}\Cc$, obtained by inverting all the quasi-isomorphisms in $\mbox{Mod-}\Cc$. This is what Keller calls the derived category of $\Cc$, and we will denote it by  $D(\modl-\Cc)$.  

There is the
standard fully faithful Yoneda embedding
\[
Z^0\Cc\to\modl-\Cc \hspace{.5in}\mbox{where}\hspace{.5in} x\in\Cc\mapsto h_x=\Cc(\cdot, x)
\]
Moreover, the Yoneda embedding induces a fully faithful functor
\[
\Ho\Cc\to D(\modl-\Cc)
\] 
This is simply because for an object $x\in\Cc$, the module $h_x$ is trivially cofibrant. 

\begin{defn}
\begin{enumerate}
\item A module $M\in \modl-\Cc$ is called {\em representable} if it is isomorphic in $\modl-\Cc$ to an object of the form $h_x$ for some $x\in\Cc$. 
\item A module $M\in \modl-\Cc$ is called {\em quasi-representable} if it is isomorphic in $D(\modl-\Cc)$ to an object of the form $h_x$ for some $x\in\Cc$. 
\end{enumerate}
\end{defn}

\begin{defn}
Let $M\in \modl-\Cc$ and $N\in \Cc-\modl$. Their tensor product is defined to be the complex 
\[
M\otimes_{\Cc} N=\mbox{cok }\{ \coprod_{c, c'\in \Cc}M(c)\otimes \Cc(c',c)\otimes N(c')\mto{\alpha} \coprod_{c\in\Cc} M(c)\otimes N(c)  \}
\]
where for $m\in M(c)$, $\phi\in \Cc(c',c)$ and $n\in N(c')$
\[
\alpha(m\otimes\phi\otimes n)=m\phi\otimes n- m\otimes \phi n
\]
\end{defn}

Bimodules are the main mechanism to construct functors between module categories over rings. They play the same role for modules over dg-categories. 
\begin{defn}
Let $\Cc$ and $\Dc$ denote two dg-categories. A bimodule $X\in \Dc-\mbox{Mod}-\Cc$ is a dg-functor 
\[
X: \Cc^\circ\otimes \Dc\to \Cc(k)
\]
More explicitly, for objects $c, c'\in\Cc$ and $d,d'\in \Dc$ there are maps of complexes 
\[
\Dc(d,d')\otimes X(c,d)\otimes\Cc(c',c)\to X(c',d')
\]
satisfying obvious conditions. 
\end{defn}
\begin{defn}
For  a bimodule $X\in \Dc-\modl-\Cc$ and $d\in \Dc$, we get an object 
\[X^d\in \modl-\Cc \hspace{.5in}\mbox{ where }\hspace{.5in} X^d(c)= X(c,d).\] 
Similarly, for $c\in \Cc$, we get an object 
\[ {^c}X\in \Dc-\modl \hspace{.5in}\mbox{ where }\hspace{.5in} {^c}X(d)= X(c,d).\]
Therefore, we may define for $M\in \modl-\Dc$ the complex
\[
M\otimes_\Dc {^c}X
\] 
Furthermore the assignment $c\mapsto {^cX}$ defines a functor $\Cc^\circ\to \Dc-\modl$ and
so $c\mapsto M \otimes_\Dc {^cX}$ defines an object in $\modl-\Cc$. Thus $\cdot \otimes_\Dc X$ defines a functor from $\modl-\Dc\to \modl-\Cc$. Moreover, by deriving this functor, we get a functor
\[
M\mapsto M\stackrel{\mathbb{L}}{\otimes}_{\Dc}X 
\]
from 
$D(\modl-\Dc)\to D(\modl-\Cc)$.
\end{defn}

\begin{defn}(Keller, \cite{Ke})
\item A bimodule $X\in \Dc-\modl-\Cc$ is called a {\em quasi-functor} if for all $d\in \Dc$, the object  
$X^d\in \modl-\Cc$ is quasi-representable.  Such a bimodule therefore defines a functor 
\[\Ho\Dc\to \Ho\Cc.\]
\end{defn}

Toen \cite{To} calls quasi-functors {\em right quasi-representable bimodules} and it is a deep theorem of his that they form the correct morphisms in the localization of the category of dg-categories by inverting dg-quasi-equivalences.

\subsection{Construction and properties of modules over $\Pc_\As$.}
We now define a class of modules over the curved dga $\As$ that will define modules over the dg-category $\Pc_\As$.

\begin{defn}
For a curved dga $\As=(\A^\bullet,d,c)$, we define a {\em quasi-cohesive module } to be the data of
$X=(X^\bullet, \mathbb{X})$ where $X^\bullet$ is a
$\Zb$-graded  right module
$X^\bullet$ over $\A$ 
together with a $\Zb$-connection
\[\mathbb{X}: X^\bullet\otimes_\A\A^\bullet\to X^\bullet\otimes_\A\A^\bullet
\]
that satisfies the integrability condition that the relative curvature 
\[F_\mathbb{X}(x)=\mathbb{X}\circ\mathbb{X}(x)+x\cdot
c=0\] for all $x\in X^\bullet$. Thus, they differ from cohesive modules  by having possibly infinitely many nonzero
graded components as well as not being projective or finitely generated over $\A$.
 
\end{defn}

\begin{defn} To a quasi-cohesive $\As$-module $X=(X^\bullet, \mathbb{X})$  we associate the $\Pc_\As$-module, denoted $\widetilde{h}_X$, by
\[
\widetilde{h}_X(E)=\{\phi:E^\bullet\otimes_\A
\A^\bullet\to X^{\bullet}\otimes_\A\A^\bullet\,\,|\,\,\mbox{ of degree $k$ and }
\phi(xa)=\phi(x)a \,\,\, \forall a\in\A^\bullet\}
\]
with differential defined in the standard way
\[
d(\phi)(ex)= \mathbb{X}(\phi(x))-(-1)^{\abs{\phi}}\phi(\Eb(x))
\]
for all $E=(E^\bullet, \mathbb{E})\in \Pc_\As$. 
We use $\widetilde{h}_X$ because of its similarity to the Yoneda embedding $h$, but beware that $X$ is not an object in $\Pc_\As$. However, in the same way as $\Pc_\As$ is shown to be a dg-category, $\widetilde{h}_X$ is shown to be a module over $\Pc_\As$. For two quasi-cohesive $\As$-modules $X$ and $Y$, and 
\[
f:X^\bullet\otimes_\A \A^\bullet \to Y^\bullet\otimes_\A \A^\bullet
\]
of degree $0$ and satisfying $f\mathbb{X}=\mathbb{Y} f$,  we get a morphism of $\Pc_\As$-modules
\[
\tilde{h}_f:\widetilde{h}_X\to \widetilde{h}_Y
\]
\end{defn}

The point of a quasi-cohesive $\As$-module  $X=(X^\bullet,\mathbb{X})$ is that the differential and morphisms decompose just the same as they do for cohesive modules. For example, $\mathbb{X}=\sum_k\mathbb{X}^k$ where $\mathbb{X}^k:E^\bullet\to X^{\bullet -k+1}\otimes_\A\A^k$ and similarly for morphisms. 
\begin{prop}
Let $X$ and  $Y$ be quasi-cohesive $\As$-modules and $f$ a morphism. Suppose $f^0:(X^\bullet, \mathbb{X}^0)\to (Y^\bullet, \mathbb{Y}^0)$ is a quasi-isomorphism of complexes. Then $\tilde{h}_f$ is a quasi-isomorphism in $\modl-\Pc_\As$. The converse is not true. 
\end{prop}  

It will be important for us to have a criterion for when a quasi-cohesive  $\As$-module $X$ induces a quasi-representable $\Pc_\As$-module. 
\begin{defn} Define a map $\phi:C\to D$ between $\A$-modules to be {\em algebraically $\A$-nuclear}, \cite{Q}, if there are finite sets of elements $\phi_k\in \Hom_\A(C,\A)$ and $y_k\in D$, $k=1,\cdots, N$ such that 
\[
\phi(x)=\sum_k y_k\cdot\phi_k(x)
\] 
\end{defn}

\begin{prop}\label{Anuc}(See Quillen, \cite{Q}, Proposition 1.1)
For $C^\bullet$ a complex of $\A$-modules, the following are equivalent:
\begin{enumerate}
\item $C^\bullet$ is homotopy equivalent to a bounded complex of finitely generated projective $\A$-modules. 
\item For any other complex of $\A$-modules $D^\bullet$, the homomorphism
\[
\Hom_\A(C^\bullet, \A)\otimes_\A D^\bullet\to \Hom_\A(C^\bullet,D^\bullet)
\]
is a homotopy equivalence of complexes (over $k$). 
\item The endomorphism $1_C$ of $C^\bullet$ is homotopic to an algebraically nuclear endomorphism.
\end{enumerate}
\end{prop}

\begin{defn}
Suppose $\As=(\A^\bullet, d, c)$ is a curved dga. Let $X=(X^\bullet,\mathbb{X})$ be a quasi-cohesive module over $\As$. Suppose there exist $\A$-linear morphisms $h^0:X^\bullet\to X^{\bullet-1}$ of degree $-1$ and $T^0:X^\bullet\to X^\bullet$ of degree $0$ satisfying
\begin{enumerate}
\item $T^0$ is algebraically $\A$-nuclear,
\item  
$
[\mathbb{X}^0,h^0]=1-T^0
$
\end{enumerate}
Then we will call $X$ a quasi-finite quasi-cohesive module. 
\end{defn}

Our criterion is the following.
\begin{thm}\label{qf}
Suppose $\As=(\A^\bullet, d, c)$ is a curved dga. Let $X=(X^\bullet,\mathbb{X})$ be a quasi-cohesive module over $\As$. 
Then there is an object $E=(E^\bullet, \mathbb{E})\in \Pc_\As$ such that  $\widetilde{h}_X$ is quasi-isomorphic to $h_{E}$; that is $\widetilde{h}_X$ is quasi-representable, under either of the two following conditions:
\begin{enumerate}
\item $X$ is a quasi-finite quasi-cohesive module.
\item $\A^\bullet$ is flat over $\A$ and there is a bounded complex $(E^\bullet,\mathbb{E}^0)$ of finitely generated projective right $\A$-modules and an $\A$-linear quasi-isomorphism $e^0:(E^\bullet,\mathbb{E}^0)\to (X^\bullet,\mathbb{X}^0)$.
\end{enumerate}
\end{thm}
\begin{proof}
In either case (1) or (2) of the theorem, there exists a bounded complex of finitely generated projective right $\A$-modules 
$(E^\bullet, \mathbb{E}^0)$ and a quasi-isomorphism $e^0:(E^\bullet,\mathbb{E}^0)\to (X^\bullet,\mathbb{X}^0)$.  In case (1), $X$ is quasi-finite-cohesive, and Proposition \ref{Anuc}, implies that $e^0$ is in fact a homotopy equivalence. In case (2) it is simply the hypothesis. 

In particular, $e^0\mathbb{E}^0-\mathbb{X}^0 e^0=0$. Now we construct a $\mathbb{Z}$-connection term by term. 
The $\mathbb{Z}$-connection $\mathbb{X}$ on $X^\bullet$ induces a connection 
\[
\mathbb{H}:H^k(X^\bullet,\mathbb{X}^0)\to H^k(X^\bullet,\mathbb{X}^0)\otimes_\A\A^1
\]
for each $k$. 
We use the quasi-morphism $e^0$ to transport this connection to a connection on $H^k(E^\bullet;\mathbb{E}^0)$
\begin{equation}
\begin{array}{ccc}
 H^k(E^\bullet;\mathbb{E}^0)& \cdots\to & H^k(E^\bullet,\mathbb{E}^0)\otimes_\A\A^1 \\
\downarrow e^0 & & \downarrow e^0\otimes 1 \\
H^k(X^\bullet,\mathbb{X}^0)& \stackrel{\mathbb{H}}\longrightarrow &  H^k(X^\bullet,\mathbb{X}^0)\otimes_\A\A^1
\end{array}\end{equation}
The right vertical arrow above $e^0\otimes 1$ is a quasi-isomorphism; in case (1) this is because $e^0$ is a homotopy equivalence and in case (2) because $\A^\bullet$ is flat. 
The first step is handled by the following lemma. 
\begin{lem}
Given a bounded complex of finitely generated projective $\A$ modules $(E^\bullet,\mathbb{E}^0)$ with  connections
$\mathbb{H}:H^k(E^\bullet;\mathbb{E}^0)\to H^k(E^\bullet,\mathbb{E}^0)\otimes_\A\A^1$, for each $k$, there exist connections
\[
\widetilde{\mathbb{H}}:E^k\to E^k\otimes_\A\A^1
\]
lifting $\mathbb{H}$. That is,  
\[
\widetilde{\mathbb{H}}\mathbb{E}^0=(\mathbb{E}^0\otimes 1)\widetilde{\mathbb{H}}
\]
and the connection induced on the cohomology is $\mathbb{H}$.
\end{lem}
\begin{proof}(of lemma) 
Since $E^\bullet$ is a bounded complex of $\A$-modules it lives in some bounded range of degrees $k\in [N,M]$. Pick an arbitrary connection on $E^M$, $\nabla$. Consider the diagram with exact rows
\begin{equation}\label{CD1.5}
\begin{array}{lclc}
  E^{M} & \mto{j} & H^M(E^\bullet,\mathbb{E}^0)&\to 0 \\
    \nabla\downarrow & \stackrel{\theta}{\searrow} &  \mathbb{H}\downarrow &\\
       & & & \\
E^{M}\otimes_\A\A^1 & \mto{j\otimes 1} & H^M(E^\bullet,\mathbb{E}^0)\otimes_\A\A^1  &\to 0
\end{array}
\end{equation}
In the diagram, $\theta=\mathbb{H}\circ j-(j\otimes 1)\circ \nabla$ is easily checked to be $\A$-linear and  $j\otimes 1$ is surjective by the right exactness of tensor product. By the projectivity of $E^M$, $\theta$ lifts to 
\[
\widetilde{\theta}:E^M\to E^M\otimes_\A\A^1
\]
so that $(j\otimes 1)\widetilde{\theta}-\theta$. 
Set $\widetilde{\mathbb{H}}=\nabla+\widetilde{\theta}$. With $\widetilde{\mathbb{H}}$ in place of $\nabla$, the diagram above commutes.

Now choose on $E^{M-1}$ any connection $\nabla_{M-1}$. But $\nabla_{M-1}$ does not necessarily satisfy $\mathbb{E}^0\nabla_{M-1}=\widetilde{\mathbb{H}}\mathbb{E}^0=0$. So we correct it as follows.  Set $\mu=\widetilde{\mathbb{H}}\mathbb{E}^0-(\mathbb{E}^0\otimes 1)\nabla_{M-1}$. Then $\mu$ is $\A$-linear. Furthermore, 
$\I \mu\subset \I \mathbb{E}^0\otimes 1$; this is because $\widetilde{\mathbb{H}}\mathbb{E}\in \I \mathbb{E}\otimes 1$ since $\widetilde{\mathbb{H}}$ lifts $\mathbb{H}$. So by  projectivity it lifts
to $\widetilde{\theta}: E^{M-1}\to E^{M-1}\otimes_\A\A^1$ such that $(\mathbb{E}^0\otimes 1)\circ\widetilde{\theta}=\theta$. Set $\widetilde{\mathbb{H}}:E^{M-1}\to E^{M-1}\otimes_\A\A^1$ to be $\nabla_{M-1}+\widetilde{\theta}$. Then $\mathbb{E}^0\widetilde{\mathbb{H}}=\widetilde{\mathbb{H}}\mathbb{E}^0$ in the right most square below. 
\begin{equation}\label{CD2}
\begin{array}{ccccccccc}
 & E^{N} & \mto{\Eb^0} & E^{N+1} &
\mto{\Eb^0}\cdots \mto{\Eb^0}& E^{M-1} & \mto{\Eb^0} & E^M&\to 0 \\
                       &   &  &  &
            & \nabla_{M-1}\downarrow & \stackrel{\mu}{\searrow} &  \widetilde{\mathbb{H}}\downarrow &\\
 & E^{N}\otimes_\A\A^1 &
\mto{\Eb^0\otimes 1} & E^{N+1}\otimes_\A\A^1 &
\mto{\Eb^0\otimes 1}\cdots\mto{\Eb^0\otimes 1} & E^{M-1}\otimes_\A\A^1 & \mto{\Eb^0\otimes 1} & E^M\otimes_\A\A^1  &\to 0
\end{array}
\end{equation}
Now we continue backwards to construct all $\widetilde{\mathbb{H}}:E^\bullet\to E^\bullet\otimes_\A\A^1$ satisfying $(\mathbb{E}^0\otimes 1)\widetilde{\mathbb{H}}=\widetilde{\mathbb{H}}\mathbb{E}^0=0$. This completes the proof of the lemma.
\end{proof}
\noindent (Proof of the theorem, continued.) Set $\widetilde{\mathbb{E}}^1=(-1)^k \widetilde{\mathbb{H}}$ on $E^k$. Then 
\[
\mathbb{E}^0\widetilde{\Eb}^1+\widetilde{\Eb}^1\mathbb{E}^0=0
\] 
but it is not necessarily true that $e^0\widetilde{\mathbb{E}}^1-\mathbb{X}^1 e^0=0$. We correct this as follows. Consider $\psi=e^0\widetilde{\mathbb{E}}^1-\mathbb{X}^1 e^0: E^\bullet\to X^\bullet\otimes_\A\A^1$. Check that $\psi$ is $\A$-linear and a map of complexes.
\begin{equation}\begin{array}{cll}
 & & (E^\bullet\otimes_\A\A^1 , \Eb^0\otimes 1)\\
 &\stackrel{\widetilde{\psi}}{\nearrow} & \downarrow e^0\otimes 1 \\
 E^\bullet & \mto{\psi} & (X^\bullet\otimes_\A\A^1, \mathbb{X}^0\otimes 1)
 \end{array}
 \end{equation}
 In the above diagram, $e^0\otimes 1$ is a quasi-isomorphism $e^0$ is a homotopy equivalence. So by Lemma 1.2.5 of \cite{OTT} there is a lift $\widetilde{\psi}$ of $\psi$ and a homotopy $e^1:E^\bullet\to X^{\bullet-1}\otimes_\A\A^1$ between $(e^0\otimes 1)\widetilde{\psi}$ and $\psi$, 
 \[
 \psi-(e^0\otimes 1)\widetilde{\psi}=(e^1\Eb^0+\mathbb{X}^0 e^1)
  \]
So let $\Eb^1=\widetilde{\Eb}^1-\widetilde{\psi}$. Then 
\begin{equation}\label{ide}
\Eb^0\Eb^1+\Eb^1\Eb^0=0 \mbox{ and } e^0\Eb^1-\mathbb{X}^1 e^0=e^1\Eb^0+\mathbb{X}^0 e^1.\end{equation}

So we have constructed the first two components $\Eb^0$ and $\Eb^1$ of the $\mathbb{Z}$-connection and the first components $e^0$ and $e^1$ of the quasi-isomorphism $E^\bullet\otimes_\A\A^\bullet \to  X^\bullet\otimes_\A\A^\bullet$. 

To construct the rest, consider the mapping cone $L^\bullet$ of $e^0$. Thus
\[
L^\bullet=E[1]^\bullet\oplus X^\bullet
\]
Let $\mathbb{L}^0$ be defined as the matrix
\begin{equation}\mathbb{L}^0=\left(\begin{array}{cc} \mathbb{E}^0[1] & 0 \\
                                         e^0[1]  & \mathbb{X}^0 \end{array}\right)
\end{equation}
Define $\mathbb{L}^1$ as the matrix
\begin{equation}\mathbb{L}^1=\left(\begin{array}{cc} \mathbb{E}^1[1] & 0 \\
                                         e^1[1]  & \mathbb{X}^1 \end{array}\right)
\end{equation}
Now $\mathbb{L}^0\mathbb{L}^0=0$ and $[\mathbb{L}^0,\mathbb{L}^1]=0$ express the identities \eqref{ide}.
Let 
\begin{equation}
D=\mathbb{L}^1\mathbb{L}^1+\left(\begin{array}{cc} 0 & 0\\
                                                    \mathbb{X}^2 e^0 & [\mathbb{X}^0,\mathbb{X}^2]
                                                    \end{array}\right) + r_c
\end{equation}
where $r_c$ denotes right multiplication by $c$. Then, as is easily checked, $D$ is $\A$-linear and  
\begin{enumerate}\item $[\mathbb{L}^0,D]=0$ and
\item $D|_{0\oplus X^\bullet}=0$.
\end{enumerate}
Since $(L^\bullet, \mathbb{L}^0)$ is the mapping cone of a quasi-isomorphism, it is acyclic and since $\A^\bullet$ is flat over $\A$, $(L^\bullet\otimes_\A\A^2, \mathbb{L}^0\otimes 1)$ is acyclic too. Since $E^\bullet$ is projective, we have that
\[
\Hom_\A^\bullet((E^\bullet,\mathbb{E}^0), (L^\bullet\otimes_\A\A^2,\mathbb{L}^0))
\]
is acyclic. Moreover \[\Hom_\A^\bullet((E^\bullet,\mathbb{E}^0), (L^\bullet\otimes_\A\A^2,\mathbb{L}^0))\subset \Hom_\A^\bullet(L^\bullet,(L^\bullet\otimes_\A\A^2,[\mathbb{L}^0,\cdot]))\] is a subcomplex. Now we have $D\in \Hom_\A^\bullet(E^\bullet, L^\bullet\otimes_\A\A^2)$ is a cycle and so there is $\widetilde{\mathbb{L}}^2\in \Hom_\A^\bullet(E^\bullet, L^\bullet\otimes_\A\A^2)$ such that $-D=[\mathbb{L}^0,\widetilde{\mathbb{L}}^2]$. Define  $\mathbb{L}^2$ on $L^\bullet$ by 
\begin{equation}
\mathbb{L}^2=\widetilde{\mathbb{L}}^2+\left(\begin{array}{cc} 0 & 0\\
                                                  0 & \mathbb{X}^2
                                                    \end{array}\right) 
\end{equation}
Then 
\begin{equation}\begin{split}
[\mathbb{L}^0,\mathbb{L}^2]= & [\mathbb{L}^0,\widetilde{\mathbb{L}}^2+\left(\begin{array}{cc} 0 & 0\\
                                                  0 & \mathbb{X}^2
                                                    \end{array}\right) ] \\
= & -D+[\mathbb{L}^0,\widetilde{\mathbb{L}}^2+\left(\begin{array}{cc} 0 & 0\\
                                                  0 & \mathbb{X}^2
                                                    \end{array}\right) ] \\
= & -\mathbb{L}^1\mathbb{L}^1-r_c
\end{split}
\end{equation}
So 
\[
\mathbb{L}^0\mathbb{L}^2+\mathbb{L}^1\mathbb{L}^1+\mathbb{L}^2\mathbb{L}^0+r_c=0.
\]

We continue by setting 
\begin{equation}
D=\mathbb{L}^1\mathbb{L}^2+\mathbb{L}^2\mathbb{L}^1+\left(\begin{array}{cc} 0 & 0\\
                                                    \mathbb{X}^3 e^0 & [\mathbb{X}^0,\mathbb{X}^3]
                                                    \end{array}\right) 
\end{equation}
Then $D:L^\bullet\to L^\bullet\otimes_\A\A^3$ is $\A$-linear, $D|_{0\oplus X^\bullet}=0$ and
\[
[\mathbb{L}^0,D]=\mathbb{L}^1\circ r_c-r_c\circ \mathbb{L}^1=0
\]
by the Bianchi identity $d(c)=0$. Hence, by the same reasoning as above, there is $\widetilde{\mathbb{L}}^3\in \Hom_\A^\bullet(E^\bullet, L^\bullet\otimes_\A\A^3)$ such that $-D=[\mathbb{L}^0,\widetilde{\mathbb{L}}^3]$. Define \begin{equation}
\mathbb{L}^3=\widetilde{\mathbb{L}}^3+\left(\begin{array}{cc} 0 & 0\\
                                                  0 & \mathbb{X}^3
                                                    \end{array}\right) 
\end{equation}
Then one can compute that $\sum_{i=0}^3\mathbb{L}^i\mathbb{L}^{3-i}=0$.

Now suppose we have defined $\mathbb{L}^0,\ldots,\mathbb{L}^n$ satisfying for $k=0,1, \ldots, n$
\[
\sum_{i=0}^k \mathbb{L}^i\mathbb{L}^{k-i}=0   \hspace{.5in}  \mbox{   for 	} k\ne 2 \\
\]
and 
\[
\sum_{i=0}^2 \mathbb{L}^i\mathbb{L}^{2-i}+r_c=0  \hspace{.5in}\mbox{   for   } k=2
\]
Then define 
\begin{equation}
D=\sum_{i=1}^n \mathbb{L}^i\mathbb{L}^{n+1-i}+\left(\begin{array}{cc} 0 & 0\\
                                                    \mathbb{X}^{n+1} e^0 & [\mathbb{X}^0,\mathbb{X}^{n+1}]
                                                    \end{array}\right) 
\end{equation}
$D|_{0\oplus X^\bullet}=0$ and we may continue the inductive construction of $\mathbb{L}$ to finally arrive at a $\mathbb{Z}$-connection satisfying $\mathbb{L}\mathbb{L}+r_c=0$. The components of $\mathbb{L}$ construct both the $\Zb$-connection on $E^\bullet$ as well as the morphism from $(E^\bullet,\mathbb{E})$ to $(X^\bullet,\mathbb{X})$. 
\end{proof}

\section{Complex manifolds}\label{Complexmanifolds}
We justify our framework in this section by showing that for a complex manifold, the derived category of sheaves on $X$ with coherent cohomology is equivalent to the homotopy category $\Pc_\As$ for the Dolbeault algebra. 
Throughout this section let $X$ be a compact complex manifold and 
$\As=(\A^\bullet, d, 0)=(\A^{0,\bullet}(X),\dbar, 0)$ the Dolbeault 
dga. This is the global sections of the sheaf of dgas $(\A_X^\bullet,d,0)=(\A_X^{0,\bullet},\dbar,0)$. Let $\mathcal{O}_X$ denote the sheaf of holomorphic functions on $X$. 
Koszul and Malgrange have shown that a
holomorphic vector bundle $\xi$ on a complex manifold $X$ is the
same thing as a $C^\infty$ vector bundle  with a flat
$\dbar$-connection, i.e., an operator
\[\dbar_\xi:E_\xi\to E_\xi\otimes_\A\A^{1}\]
such that $\dbar_\xi(f\phi)=\dbar(f)\phi+f\dbar_\xi(\phi)$ for $f\in \A, \phi\in \Gamma(X;\xi)$
and satisfying the integrability condition that $\dbar_\xi\circ\dbar_\xi=0$. Here $E_\xi$ denotes the global $C^\infty$ sections of $\xi$. The notion of a cohesive module over $\As$ clearly generalizes this notion but in fact will also include coherent analytic
sheaves on $X$ and even more generally, bounded complexes of $\mathcal{O}_X$-modules with coherent cohomology as well. 

For example, if  $(\xi^\bullet, \delta)$ denotes a complex of holomorphic vector bundles, with corresponding global $C^\infty$-sections $E^\bullet$ and $\dbar$-operator $\dbar_\xi:E^i\to E^i\otimes_\A\A^1$ then the holomorphic condition on $\delta$ is that $\delta\dbar_\xi=\dbar_\xi\delta$. Thus $E=(E^\bullet,\mathbb{E})$, where $\mathbb{E}^0=\delta$ and $\mathbb{E}^1=(-1)^\bullet\dbar_\xi$ defines the cohesive module corresponding to $(\xi^\bullet,\delta)$. So we see that for coherent sheaves with locally free resolutions, there is nothing new here.

\subsection{The derived category of sheaves of $\mathcal{O}_X$-modules with coherent cohomology}
Pali, \cite{Pa} was the first to give a characterization of general coherent analytic sheaves in terms of sheaves over $(\A_X^\bullet, d)$ equipped with flat $\dbar$-connections. 
He defines a {\em $\dbar$-coherent
analytic sheaf} $\mathcal{F}$ to be a sheaf of modules over the
sheaf of $C^\infty$-functions ${\A_X}$ satisfying two
conditions:
\begin{enumerate}
\item Finiteness: locally on $X$, $\mathcal{F}$ has a finite
resolution by finitely generated free modules, and \item Holomorphic: $\mathcal{F}$ is equipped with a $\dbar$-connection, i.e., an operator (at the level of sheaves)
\[\dbar:\mathcal{F}\to \mathcal{F}\otimes_{{\A}_X}{\A}_X^{1}\]
and satisfying $\dbar^2=0$.
\end{enumerate}
\begin{thm}(Pali, \cite{Pa})\label{Pa}
The category of coherent analytic sheaves on $X$ is equivalent to
the category of $\dbar$-coherent sheaves.
\end{thm}

We prove our theorem independently of his. We use the following proposition of Illusie, \cite{SGA6}.
\begin{prop}\label{SoftStuff}
Suppose $(X,\A_X)$ is a ringed space, where $X$ is compact and $\A_X$ is a soft sheaf of rings. Then 
\begin{enumerate}
\item The global sections functor 
\[
\Gamma:\text{Mod-}\A_X\to \text{Mod-}\A_X(X)
\]
is exact and
establishes an equivalence of categories between the category of sheaves of right $\A_X$-modules and the category of right modules over the global sections $\A_X(X)$.
\item If $M\in\text{Mod-}\A_X$ locally has finite resolutions by finitely generated free $\A_X$-modules, then $\Gamma(X;M)$ has a finite resolution by finitely generated projectives. 
\item The derived category of perfect complexes of sheaves $D_{\mbox{perf}}(\text{Mod-}\A_X)$ is equivalent 
the derived category of perfect complexes of modules $D_{\mbox{perf}}(\text{Mod-}\A_X(X))$.
\end{enumerate}
\end{prop}
\begin{proof} See Proposition 2.3.2, Expos\'{e} II, SGA6, \cite{SGA6}. \end{proof}

Our goal is to derive the following description of the bounded derived category of sheaves of $\mathcal{O}_X$-modules with coherent cohomology on a complex manifold. Note that this is equivalent to the category of perfect complexes, since we are on a smooth manifold. Recall that $\As=(\A^\bullet, d, 0)=(\A^{0,\bullet}(X),\dbar, 0)$ the Dolbeault 
dga is the global sections of the sheaf of dgas $(\A_X^\bullet,d,0)=(\A_X^{0,\bullet},\dbar,0)$

\begin{thm}\label{cateq} Let $X$ be a compact complex manifold and $\As=(\A^\bullet, d, 0)=(\A^{0,\bullet}(X),\dbar, 0)$ the Dolbeault 
dga. Then the category $\Ho\Pc_{\As}$ is equivalent
to the bounded derived category of complexes of sheaves of $\mathcal{O}_X$-modules with
coherent cohomology ${D}^b_{\mbox{coh}}(X)$.
\end{thm}
\begin{rem} This theorem is stated only for $X$ compact. This is because Proposition \ref{SoftStuff} is stated only for $X$ compact. A version of Theorem \ref{cateq} will be true once one is able to characterize the perfect $\A_X$-modules in terms of modules over the global sections for $X$ which are not compact. 
\end{rem}

A module $M$ over $\A$ naturally localizes to a sheaf $M_X$ of $\A_X$-modules where 
\[
M_X(U)=M\otimes_\A \A_X(U)
\]

For an object $E=(E^\bullet,\Eb)$ of $\Pc_{\As}$, define the
sheaves $\mathcal{E}_X^{p,q}$ by \[\mathcal{E}_X^{p,q}(U)=
E^p\otimes_{\A}\A^q_X(U).\] We define a complex of sheaves by
$(\mathcal{E}_X^\bullet,\mathbb{E})=(\sum_{p+q=\bullet}\mathcal{E}_X^{p,q},\Eb)$.
This is a complex of soft sheaves of $\mathcal{O}_X$-modules, since $\Eb$ is a $\dbar$-connection. 
The theorem above will be broken up into several lemmas.

\begin{lem}\label{lem1}
The complex $\mathcal{E}_X^\bullet$ has coherent cohomology and 
\[
E=(E^\bullet,\Eb)\mapsto \alpha(E)=(\mathcal{E}_X^\bullet,\Eb)
\] 
defines a fully faithful functor $\alpha:\Ho\Pc_{\As}\to \mathcal{D}_{\mbox{perf}}(X)\,\widetilde{-}\,\mathcal{D}^b_{\mbox{coh}}(X)$.
\end{lem}

\begin{proof} Let $U$ be a polydisc in $X$. We show that on a possibly smaller polydisc $V$, there is gauge transformation 
$\phi:\mathcal{E}^\bullet|_V\to\mathcal{E}^\bullet|_V$ of degree zero such that $\phi\circ\Eb\circ \phi^{-1}=\mathbb{F}^{0}+\dbar$. Thus $\mathcal{E}^\bullet|_V$ is gauge equivalent to a complex of holomorphic vector bundles.  Or in other words, for each $p$
the sheaf $H^p((\mathcal{E}^{ \bullet,0},\Eb^0)$ is
$\dbar$-coherent, with $\dbar$-connection $\Eb^1$. Since $U$ is
Stein there is no higher cohomology (with respect to $\Eb^1$) and
we are left with the holomorphic sections over $U$ of each of
these $\dbar$-coherent sheaves, which are thus coherent. 

The construction of the gauge transformation follows the proof of the integrability theorem for complex structures on vector bundles, \cite{DK}, section 2.2.2, page 50.
Thus we may assume we are in a polydisc $U=\{(z_1,\cdots,z_n)|\,\, |z|_i < r_i\}$. In these coordinates we may write the $\mathbb{Z}$-connection $\Eb$ as $\Eb=\Eb^0+\dbar+J$ where 
\[
J:\mathcal{E}^{p,q}(U)\to\bigoplus_{i\le p}\mathcal{E}^{i,q+(p-i)+1}(U)
\]
is $\A_X(U)$-linear. Now write $J=J'\wedge d\bar{z}_1+J''$ where $\iota_{\frac{\partial}{\partial \bar{z}_1}}J'=\iota_{\frac{\partial}{\partial \bar{z}_1}}J''=0$. Write $\dbar_i$ for $d\bar{z}_i\wedge\frac{\partial}{\partial \bar{z}_i}$. As in \cite{DK}, page 51, we find a $\phi_1$ such that 
$\phi_1(\dbar_1+J'\wedge d\bar{z}_1)\phi_1^{-1}=\dbar_1$, by solving $\phi_1^{-1}\dbar_1(\phi_1)=J'\wedge d\bar{z}_1$, for $\phi_1$, possibly having to shrink the polydisc.  Here, we are treating the variables $z_2,\cdots z_n$ as parameters. Then we set $\Eb_1=\phi_1(\Eb^0+\dbar +J'+J'')\phi_1^{-1}$. Then $\Eb_1\circ\Eb_1=0$ and we can write 
\[
\Eb_1=\Eb_1^0+\dbar_1+\dbar_{\ge 2}+J_1
\]
where $\iota_{\frac{\partial}{\partial \bar{z}_1}}J_1=0$ and we can check that both $\Eb_1^0$ and $J_1$ are holomorphic in $z_1$. For $0=\Eb_1\circ\Eb_1$ and therefore
\begin{equation}\begin{split}
0&=\iota_{\frac{\partial}{\partial \bar{z}_1}}(\Eb_1\circ\Eb_1)\\
  & = \iota_{\frac{\partial}{\partial \bar{z}_1}}(\Eb_1^0\circ\dbar_1+\dbar_1\circ E_1^0+J_1\circ\dbar_1+\dbar_1\circ J_1 )\\
   &=\iota_{\frac{\partial}{\partial \bar{z}_1}}(\dbar_1(\Eb_1^0)+\dbar_1(J_1))
  \end{split}\end{equation}
Now each of the two summands in the last line must individually be zero since $\iota_{\frac{\partial}{\partial \bar{z}_1}}(\dbar_1(\Eb_1^0))$ increases the $p$-degree by one and $\iota_{\frac{\partial}{\partial \bar{z}_1}}(\dbar_1(J_1))$ preserves or decreases the $p$-degree by one. So we have arrived at the following situation:
\begin{enumerate}
\item $\Eb_1^0\circ\Eb_1^0=0$,
\item $\Eb_1^0$ and $J_1$ are holomorphic in $z_1$, and 
\item $\iota_{\frac{\partial}{\partial \bar{z}_1}}J_1=0$.
\end{enumerate}
We now iterate this procedure. Write $J_1=J_1'\wedge d\bar{z}_2+J_1''$ where $\iota_{\frac{\partial}{\partial \bar{z}_1}}J_1'=\iota_{\frac{\partial}{\partial \bar{z}_2}}J_1'=\iota_{\frac{\partial}{\partial \bar{z}_1}}J_1''=\iota_{\frac{\partial}{\partial \bar{z}_2}}J_1''=0$. Now solve 
\[
\phi_2^{-1}\dbar_2(\phi_2)=J_1'\wedge d\bar{z}_2
\]
for $\phi_2$. Since $J_1'$ is holomorphic in $z_1$ and smooth in $z_2,\cdots, z_n$, so will $\phi_2$. Then as before we  have 
\[
\phi_2(\dbar_2+J_1'\wedge d\bar{z}_2)\phi_2^{-1}=\dbar_2
\]
as well as 
\[
\phi_2(\dbar_1)\phi_2^{-1}=\dbar_1
\]
since $\phi_2$ is holomorphic in $z_1$.
Setting $\Eb_2=\phi_2\circ\Eb_1\circ\phi_2^{-1}$, we see that 
\[
\Eb_2=\Eb_2^0+\dbar_1+\dbar_2+\dbar_{\ge 3}+J_2
\]
where $\iota_{\frac{\partial}{\partial \bar{z}_1}}J_2=\iota_{\frac{\partial}{\partial \bar{z}_2}}J_2=0$ and we can check as before that both $\Eb_2^0$ and $J_2$ are holomorphic in $z_1$ and $z_2$. And continue until we arrive at $\mathbb{F}=\Eb_{n}=\Eb_n^0+\dbar$.

  \end{proof}

\begin{lem}\label{lem2}
To any complex of sheaves of $\mathcal{O}_X$-modules $(\Ec_X^\bullet, d)$ on $X$ with coherent cohomology, there
corresponds a cohesive $\As$-module $E=(E^\bullet, \mathbb{E})$, unique up to
quasi-isomorphism in $\Pc_{\As}$ and a quasi-isomorphism 
\[
\alpha(E)\to (\Ec^\bullet, d)
\]
This correspondence has
the property that for any two such complexes 
$\Ec_1^\bullet$ and $\Ec_2^\bullet$, that the corresponding
twisted complexes $(E_1^\bullet,\Eb_1)$ and $(E_2^\bullet,\Eb_2)$
satisfy
\[
\Ext^k_{\mathcal{O}_X}(\Ec_1^\bullet,\Ec_2^\bullet)\cong
H^k(\Pc_\As(E_1,E_2))
\]
\end{lem}
\begin{proof}
Since we are on a manifold we may assume that $(\Ec^\bullet,d)$ is a perfect complex. Set $\Ec^\bullet_\infty =
\Ec^\bullet\otimes_{\mathcal{O}_X}\A_X$. Now the map $(\Ec^\bullet,d)\to (\Ec^\bullet_\infty\otimes_\A\A^{\bullet}_X,d\otimes 1+1\otimes \dbar)$ is a quasiisomorphism of sheaves of $\mathcal{O}_X$-modules by the flatness of $\A_X$ over $\mathcal{O}_X$. Again, by the flatness of $\A_X$ over $\mathcal{O}_X$, it follows that $(\Ec^\bullet_\infty,d)$ is a perfect complex of $\A_X$-modules. By Proposition \ref{SoftStuff}, there is a (strictly) perfect complex $(E^\bullet, \mathbb{E}^0)$ of $\A$-modules and quasiisomorphism $e^0:(E^\bullet,\mathbb{E}^0)\to (\Gamma(X,\Ec^\bullet_\infty),d)$. Moreover $(\Gamma(X,\Ec^\bullet_\infty),d\otimes1+1\otimes \dbar)$ defines a quasi-cohesive module over $\As$. So the hypotheses of Theorem \ref{qf}(2) are satisfied. The lemma is proved.  \end{proof}

\subsection{Gerbes on complex manifolds}\label{gerbedescription}
The theorem above has an analogue for gerbes over compact manifolds. $X$ is still  compact complex manifold. A class $b\in H^2(X,\mathcal{O}_X^\times)$ 
defines an $\mathcal{O}_X^\times$-gerbe on $X$. From the exponential sequence of sheaves 
\[0\to \Zb_X\to \mathcal{O}_X\mto{\exp \cf\cdot} \mathcal{O}_X^\times\to 0\]
there is a long exact sequence
\[
\cdots\to H^2(X;\mathcal{O}_X)\to H^2(X;\mathcal{O}_X^\times)\to H^3(X;\Zb_X)\to\cdots 
\]
If $b$ maps to $0\in H^3(X;\Zb)$ (that is, the gerbe is topologically trivializable) then $b$ pulls back to a class represented by a $(0,2)$-form $B\in \A^{0,2}(X)$. Consider the curved dga $\As=(\A^\bullet, d, B)=(\A^{0,\bullet}(X),\dbar,B)$; the same Dolbeault algebra as before but with a curvature. 
Then we have a theorem \cite{BD}, corresponding to \ref{cateq}, 
\begin{thm}\label{cateq2} The category $\Ho\Pc_{\As}$ is equivalent
to the bounded derived category of complexes of sheaves on the gerbe $b$ over $X$ of $\mathcal{O}_X$-modules with
coherent cohomology and weight one ${D}^b_{{coh}}(X)_{(1)}$.
\end{thm}
Sheaves on a gerbe are often called twisted sheaves. One can deal with gerbes which are not necessarily topologically trivial, but the curved dga is slightly more complicated, \cite{BD}.

\section{Examples}
\subsection{Elliptic curved dgas}
In this section we define a class of curved dga's $\As$ such that the corresponding dg-category $\Pc_\As$ is {\em proper}, that is, the cohomology of the hom sets are finite dimensional. It is often useful to equip a manifold with a Riemannian metric so that one can use Hilbert space methods. We introduce a relative of the notion spectral triple in the sense of Connes, \cite{Co1}, so that we can use Hilbert space methods to guarantee the properness of the dg-category. 

Again, our basic data is a curved dga $\As=(\A^\bullet, d,c)$. 
\begin{defn}
We say that $\As$ is equipped with a {\em Hilbert structure} if there is a positive definite Hermitian inner product on $\A^\bullet$
\[
\langle \cdot,\cdot\rangle:\A^k\times\A^k\to\mathbb{C}
\]
satisfying the following conditions: Let $\Hc^\bullet$ be the completion of $\A^\bullet$. 
\begin{enumerate}
\item For $a\in \A^\bullet$, the operator $l_a$ (respectively  $r_a$) of left (respectively, right) multiplication by $a$ extends to $\Hc^\bullet$ as a bounded operator. Furthermore, the operators $l_a^*$ and $r_a^*$ map $\A^\bullet\subset \Hc^\bullet$ to itself. 
\item $\A$ has an anti-linear involution $*:\A\to \A$ such that for $a\in\A$, there is $(l_a)^*=l_{a^*}$ and $(r_a)^*=r_{a^*}$.
\item The differential $d$ is required to be closable in $\Hc^\bullet$. Its adjoint satisfies  $d^*(\A^\bullet)\subset \A^\bullet $ and the operator $D=d+d^*$ is essentially self-adjoint with core $\A^\bullet$.
\item For $a\in\A^\bullet$, $[D,l_a]$, $[D,r_a]$, $[D,l_a^*]$ and $[D,r_a^*]$ are bounded operators on $\Hc^\bullet$.
\end{enumerate}
\end{defn}

\begin{defn} An elliptic curved dga $\As=(\A^\bullet, d,c)$ is a curved dga with a Hilbert structure which in addition satisfies
\begin{enumerate}
\item   The operator $e^{-tD^2}$ is trace class for all $t>0$. 
\item \[\A^\bullet = \bigcap_n \mbox{Dom}(D^n)\] 
\end{enumerate}
\end{defn}

\noindent The following proposition follows from very standard arguments. 
\begin{prop} Given an elliptic curved dga $\As$ then 
for $E=(E^\bullet,\Eb)$ and $F=(F^\bullet,\Fb)$ in $\Pc_\As$ one has that the cohomology of $\Pc_\As(E,F)$ is finite dimensional. 
\end{prop}

Bondal and
Kapranov have given a very beautiful formulation of Serre duality
purely in the derived category. We adapt their definitions to our
situation of dg-categories.
\begin{defn}
For a dg-category $\Cc$, such that all $\Hom$ complexes have
finite dimensional cohomology, a Serre functor is a dg-functor
\[
S:\Cc\to \Cc
\]
which is a dg-equivalence and so that there are 
pairings of degree zero, functorial in both $E$ and $F$
\[
\langle\cdot,\cdot\rangle:\Cc^\bullet(E,F)\times
\Cc^\bullet(F,SE)\to \mathbb{C}[0]
\]
satisfying 
\[
\langle d\phi,\psi\rangle+(-1)^{|\phi|}\langle \phi,d\psi\rangle =
0
\]
which are perfect on cohomology for any $E$ and $F$ in $\Cc$.
\end{defn}

Motivated by the case of Lie algebroids below, we make the following definition, which will guarantee the existence of a Serre functor.

\begin{defn}
Let $\As=(\A^\bullet, d, c)$ be an elliptic curved dga. A dualizing module (of dimension $g$) is a triple $((D,\mathbb{D}), \bar{*}, \int)$ where
\begin{enumerate}
\item $(D, \mathbb{D})$ is an $\As-\As$ cohesive  bimodule,
\item $\bar{*}:\A^k\to D\otimes_\A\A^{g-k}$ is a conjugate linear
isomorphism and satisfies
\[
\bar{*}(a\omega)=\bar{*}(\omega)a^* \hspace{.5in} \mbox{and} \hspace{.5in} \bar{*}(\omega a)=a^*\bar{*}(\omega)
\]
for $a\in\A$ and $\omega\in \A^\bullet.$
 \item a $\mathbb{C}$-linear map
$\int:D\otimes_\A\A^g\to \mathbb{C}$ such that $\int
\mathbb{D}(x)=0$ for all $x\in D\otimes_\A\A^\bullet$ and
\[
\int \omega\cdot x=(-1)^{|\omega||x|}\int x\cdot\omega
\]
for all $\omega\in \A^\bullet$ and $x\in D\otimes_\A\A^\bullet$,
and \[\langle \omega,\eta\rangle = \int \bar{*}(\omega)\eta
\]
\end{enumerate}
\end{defn}

\begin{prop}
Given an elliptic curved dga $\As=(\A^\bullet, d, c)$ with a dualizing module $((D,\mathbb{D}), \bar{*},
\int)$, the category $\Pc_{\As}$ has a Serre functor given
by the cohesive bimodule $(D[g],\mathbb{D})$. That is,
\[
S(E^\bullet,\mathbb{E})= (E\otimes_\A D[g],
\mathbb{E}\#\mathbb{D})
\]
is a dg-equivalence for which there are functorial pairings
\[
\langle\cdot,\cdot\rangle:\Pc_{\As}^\bullet(E,F)\times
\Pc_{\As}^\bullet(F,SE)\to \mathbb{C}
\]
satisfying
\[
\langle d\phi,\psi\rangle+(-1)^{|\phi|}\langle \phi,d\psi\rangle =
0
\]
is perfect on cohomology for any $E$ and $F$ in
$\Pc_{\As}$.
\end{prop}

\subsection{Lie Algebroids}

Lie algebroids provide a natural source of dga's and
thus, by passing to their cohesive modules, interesting dg-categories. 

Let $X$ be a $C^\infty$-manifold and let $\Al$ be a complex Lie algebroid over $X$. Thus $\Al$ is a $C^\infty$ vector
bundle on $X$ with a bracket operation on $\Gamma(X;\Al)$
making $\Gamma(X;\Al)$ into a Lie algebra and such that the induced map into
vector fields $\rho:\Gamma(X;\Al)\to \V(X)$ is a Lie algebra homomorphism and
for $f\in C^\infty(X)$ and $x,y\in \Gamma(X;\Al)$ we have
$$ [x,fy]=f[x,y]+(\rho(x)f)y$$
Let $g$ be the rank of $\Al$ and $n$ for the
dimension of $X$.

There is a dga corresponding to any Lie algebroid $\Al$ over $X$
as follows. Let
$$\A_\Al^\bullet =\Gamma(X;\wedge^\bullet \Al^\vee)$$ denote the space of
smooth $\Al$-differential forms. It has a differential $d$ of degree
one, with $d=0$ given by the usual formula,
\begin{eqnarray}
\label{eq_dA}
 (d \eta)  (x_1, ..., x_k) & = &  \sum_i
 (-1)^{i+1} \rho(x_i) (\eta (x_1, ..., \hat{x}_i, ..., x_k)) \\ &   + &
 \sum_{i < j} (-1)^{i+j} \eta([x_i, x_j], ..., \hat{x}_i, ...,
 \hat{x}_j, ..., x_k). \nonumber
 \end{eqnarray}
  turning it into a
differential graded algebra. Note that $\A_\Al=\A_\Al^0$ is
just the $C^\infty$-functions on $X$. Then $\As_\Al=(\A_\Al^\bullet,d,0)$ is a curved dga.

 \subsubsection{The dualizing  $\Al$-module $\D_\Al$}

 \bigskip\noindent
 We recall the definition of the ``dualizing module" of a Lie
 algebroid. This was first defined in \cite{ELW} where they used it to
 define the modular class of the Lie algebroid.

 \bigskip\noindent
Let $\Al$ be a Lie algebroid over $X$ with anchor map $\rho$.  Consider
 the line bundle
 \begin{equation}
 \label{eq_qa}
 \D_{\Al} ~ = ~  \wedge^g \Al ~ \otimes ~  \wedge^n
 T^\vee_\cx X.
 \end{equation}
 Write $D_\Al=\Gamma(X;\D_\Al)$.
 Define
 \[
 \mathbb{D}: ~~ D_\Al\to D_\Al\otimes_{\A_\Al}\A_\Al^1
 \]
 by
 \begin{equation}
 \label{eq_d-on-qa}
 \mathbb{D}(X \otimes \mu)(x) ~ = ~ L_x (X) ~ \otimes ~ \mu ~ + ~ X ~ \otimes ~
 L_{\rho(x)} \mu,
 \end{equation}
 where $x\in \Gamma(X;\Al)$, $X \in \Gamma(X;\wedge^g\Al), ~  \mu \in  \Gamma(\wedge^n
 T^\vee_\cx X)$, and $L_{\rho(x)}\mu$ denotes the Lie derivative of $\mu$ in
 the direction of $\rho(x)$. See \cite{ELW} for more details.

Now we note that $\As_\Al$ acts on the left of $D_\Al\otimes_{\A_\Al}\A_\Al^\bullet$ and
$\mathbb{D}:D_\Al\to D_\Al\otimes_{\A_\Al}\A_\Al^1$
defines a flat $\A_\Al^\bullet$-connection \cite{ELW}.
Therefore $(D_\Al,\mathbb{D})$ denote a cohesive $\As_\Al-\As_\Al$-bimodule, and thus a dg-functor from  $\Pc_{\As_\Al}$ to itself. 

\noindent We have the pairing
$$\D_\Al\otimes  \wedge^{g}\Al^\vee\to\wedge^n T_\cx^*X.$$ Which allows us to define 
 $\int:D_\Al\otimes_\A\A^g\to \mathbb{C}$ for  $(X\otimes\mu)\otimes\nu\in
  \D_\Al\otimes\wedge^{g}\Al^\vee $
 \[
 \int (X,\nu)\mu
 \]
 Then we have
 \begin{thm}[Stokes' Theorem, \cite{ELW}]
 \label{thm_stokes}
 Identify $D_\Al\otimes_{\A_\Al}\A_\Al^g(X) = \Gamma(\wedge^g \Al \otimes \wedge^g \Al^\vee \otimes
 \wedge^nT_\cx^\vee X)$ with the space of top-degree forms on $X$ by
 pairing the factors in $\wedge^g \Al^*$ and $\wedge^g \Al$ pointwise. We
 have, for  every $c = (X\otimes\mu)\otimes\nu\in D_{\A_\Al}\otimes_{\A_\Al}\A_\Al^{g-1}(X)$,
 \begin{equation}
 \label{eq_a-a-r-1}
 \mathbb{D}(c) ~ = ~ (-1)^{g-1} d( \rho(\mu \lrcorner X) \lrcorner \nu).
 \end{equation}
 Consequently,
 \begin{equation}
 \label{eq_int-0}
 \int_X \mathbb{D}(c) ~ = ~ 0.
 \end{equation}
 \end{thm}

\subsubsection{Hermitian structures and the $\bar{*}$-operator}
 Let us equip the  algebroid $\Al$ with an Hermitian inner
 product $< ~,~>$. Then $\Al^\vee$ and $\wedge^\bullet \Al^\vee$ all
 inherit Hermitian inner products according to the rule
$$\langle\alpha_1\wedge\cdots\wedge\alpha_k, \,\,
 \beta_1\wedge\cdots\wedge\beta_k\rangle=\mbox{det} (\langle\alpha_i,\beta_j\rangle).$$
 Also let us put on $X$ a Riemannian structure and let $\nu_X$ be
 the volume form. Then there is a Hermitian inner product on
 $\A_\Al^\bullet(X)$ defined by
 $$\langle\alpha,\,\,\beta\rangle=\int_X \langle\alpha,\,\,\beta\rangle\nu_X.$$
 Recall that there is a canonical identification of $ \D_\Al\otimes \wedge^g\Al^\vee$ with
 $\Lambda^n T_\cx^\vee X$.
 Define the operator $\bar{*}:\wedge^k\Al^\vee\to
 \D_\Al\otimes  \wedge^{g-k}\Al^\vee$ by requiring that
 \begin{equation}\label{stardefn}\alpha\wedge \bar{*}\beta=\langle\alpha,\,\,
 \beta\rangle\nu_X.\end{equation}
This is well defined because the pairing
$$\wedge^k\Al^\vee\times  (\D_\Al\otimes \wedge^{g-k}\Al^\vee)\to \wedge^nT_\cx^\vee X$$ is perfect.
Our $\bar{*}$ operator is conjugate linear. This is because we
have no conjugation operator on $\Al$, as would be the case when
we define the Hodge $*$ operator on the bigraded Dolbeault
complex.

\bigskip\noindent As usual, we have the familiar local expressions for the
$\bar{*}$-operator. So if $\alpha_1,\cdots,\alpha_g$ is an
orthonormal frame of $\Al$ with $\alpha^1,\cdots,\alpha^g$ the
dual frame, then for a multi-index $I\subset \{1,\ldots,g \}$ we
have
$$\bar{*}(\lambda\alpha^I)=(-1)^{\sigma(  I ) }
(\alpha_{\{1,\ldots,g\}}\otimes\nu_X)\otimes
\bar{\lambda}\alpha^{I^c}$$ where $I^c$ is the complement of the
multiindex and $\sigma(I)$ is the sign of the permutation
$(1,\ldots,g)\mapsto (I,I^c)$. For an object $(E,\mathbb{E})$ of
$\Pc_{\As_\Al}$ we equip
$E$ with a Hermitian structure (no condition) and we extend
$\bar{*}$ to
$$\bar{*}: E\otimes\wedge^k\Al^\vee\to
 E^\vee\otimes \D_\Al\otimes\wedge^{g-k}\Al^\vee $$ by the same
 formula, (\ref{stardefn}). Locally, we have
 $$\bar{*}( e_i\otimes \lambda\alpha^I)=(-1)^{\sigma(  I ) }
 e^i\otimes(\alpha_{\{1,\ldots,g\}}\otimes
 \nu_X)\otimes \bar{\lambda}\alpha^{I^c}
 $$
 where $e_i$ and $e^i$ are dual pairs of orthonormal frames of $E$
 and $E^\vee$ respectively.

\np Now we make a basic assumption on our Lie algebroid $\Al$.
\begin{defn}
A complex Lie algebroid $\rho:\Al\to T_\cx X$ is called {\em
elliptic} if $$\rho^\vee:T^\vee X\to T_\cx^\vee X\to\Al^\vee$$ is
injective.
\end{defn}

\np Note that a real Lie algebroid is elliptic means that it is
transitive.  The point of this definition is the following
proposition.
\begin{thm}
For an elliptic Lie algebroid $\Al$, the corresponding dga $\As_\Al=(\A_\Al^\bullet, d, 0)$ is an elliptic dga and $((D_\Al,\mathbb{D}), \bar{*},\int)$ is a dualizing manifold. with a representation
\end{thm}
\begin{proof}
Everything follows from basic elliptic theory. 
\end{proof}

\np As an immediate corollary we have
\begin{thm}\label{dualitytheorem}
For an elliptic Lie algebroid $\Al$ with $(E,\mathbb{E}), (F,\mathbb{F})\in \Pc_{\As_\Al}$, there is a perfect duality pairing 
$$H^k(\Pc_{\As_\Al}(E, F))\times H^{g-k}(\Pc_{\As_\Al}(F,E\otimes D_\Al))\to \cx  $$
\end{thm}

\subsubsection{The De Rham Lie algebroid and Poincar\'{e} duality}\label{deRham} For
$\rho=\mbox{Id}:\Al=TM\to TM$ the duality theorem is
Poincar\'{e}'s for local systems. That is, the dualizing module is
the trivial one dimensional vector bundle (we made the blanket
assumption that $M$ is orientable) and for a flat vector bundle
$E$ over $X$ there is a perfect pairing
$$H^k(X;E)\times H^{n-k}(X,E^\vee)\to \cx$$

\subsubsection{The Dolbeault Lie algebroid and Serre
Duality}\label{Dolbeault}
 For $X$ a complex $n$-dimensional manifold, let
$\rho:\Al=T^{0,1}\hookrightarrow T_\cx X$ be the natural
inclusion. Thus, a holomorphic vector bundle is the same thing as
an $T^{0,1}$-module. Moreover
$$\D_{T^{0,1}} =   \wedge^n T^{0,1} ~ \otimes  \wedge^{2n}
T^\vee_\cx X\cong \wedge^n T^{0,1} ~ \otimes \wedge^n
(T^{0,1\,\,\vee}X\oplus T^{1,0\,\,\vee } X)$$
$$\cong \wedge^n T^{1,0\,\,\vee}X$$
is the usual canonical (or dualizing) bundle $K$ in complex
geometry. And  (\ref{dualitytheorem}) reduces to Serre's duality
theorem that for a holomorphic vector bundle $E$ the sheaf (i.e.
Dolbeault) cohomology satisfies
$$H^{k}_{\bar{\partial}}(X;E)^\vee\cong H^{n-k}_{\bar{\partial}}(X;E^\vee\otimes K).$$
from which it follows by letting $E$ be $\wedge^p {T^{1,0}}^\vee$
that
$$H^{p,q}(X)\cong H^{n-p,n-q}(X).$$
Stated in terms of Serre functors, we have that with $SE=E\otimes
K[n]$, $S$ is a Serre functor on $\Pc_{T^{0,1}}$.

\subsubsection{The Higgs Lie algebroid}\label{Higgs}
Again, let $X$ be an $n$-dimensional complex manifold. We define a
new Lie algebroid as follows.
$$\Al=T_\cx X=T^{0,1}\oplus
T^{1,0}\stackrel{p''}{\rightarrow} T_\cx X$$
 where $p''$ is the
projection of the complexified tangent bundle onto $T^{0,1}X$. Let
$p'$ be the projection onto $T^{1,0}X$. We need to adjust the
bracket by
\[
\{X'+X'',Y'+Y''\}=[X'',Y'']+p'([X'',Y']+[X',Y'']).
\]
for $X', Y'\in\Gamma (T^{1,0})$ and $X'', Y''\in \Gamma (T^{0,1})$
and the square brackets denote the usual bracket if vector fields.
\begin{prop}
\begin{enumerate}

\item $\Al$ is an elliptic Lie algebroid.

\item A module over $\Al$ is comprised of the following data:
$(E,\Phi)$ where $E$ is a holomorphic vector bundle and $\Phi$ is
a holomorphic section of $\mbox{Hom}(E,E\otimes T^{1,0\vee}X)$ and
satisfies the integrability condition $\Phi\wedge\Phi=0$, that is
$(E,\Phi)$ is a Higgs bundle in the sense of Hitchin, and Simpson, \cite{Hit2}, 
\cite{Sim1}.

\item The dualizing module $\D_\Al$ is the trivial one dimensional
vector bundle with the Higgs field $\Phi=0$. \end{enumerate}
\end{prop}
\proof That $\{\cdot,\cdot\}$ satisfies Jacobi is a
straightforward calculation that only uses the integrability of
the complex structure, that is, that $T^{0,1}$ and $T^{1,0}$ are
both closed under bracket. To check that algebroid condition we
calculate $\{X,fY\}$ \begin{equation}\begin{array}{l}
 =\{X'+X'',fY'+fY''\} = [X''+fY'']+p'([X',fY'']+[X'',fY']) \\
=f[X'',Y'']+X''(f)Y''+p'(f[X',Y''] +X'(f)Y''+f[X'',Y']+X''(f)Y')
\\
=f([X'',Y'']+p'([X',Y'']+[X'',Y']))
+X''(f)(Y'+Y'')+p'(X'(f)Y'')\\
=f\{X,Y\}+p''(X)(f)Y
\end{array}\end{equation}
To show it is an elliptic Lie algebroid, let $\xi\in TX^\vee$.
Since its image in $T^\vee_\cx X$ is real it can be written as
$e+\bar{e}$ for $e\in T^{0,1\vee}$. The projection to
$T^{0,1\vee}$ is simply $e$ and thus $\rho^\vee$ is injective from
$TX^\vee\to \Al^\vee$.

For the statement about modules, suppose $(E,\mathbb{E})$ is a
module over $\Al$. Then we have the decomposition
\[
\mathbb{E}:\Gamma(E)\to \Gamma(E\otimes (T^{1,0\vee}\oplus
T^{0,1\vee}))\cong  \Gamma(E\otimes T^{1,0\vee})\oplus
\Gamma(E\otimes T^{0,1\vee})
\]
in which $\mathbb{E}$ decomposes as $\mathbb{E}=\mathbb{E}'\oplus
\mathbb{E}''$. The condition of being an $\Al$ connection means
that $\mathbb{E}''$ satisfies Leibniz with respect to the
$\dbar$-operator and $\mathbb{E}'$ is linear over the functions
and thus $\Phi=\mathbb{E}':E\to E\otimes T^{1,0\vee}$. The
flatness condition $\nabla^2=0$ implies
\begin{enumerate}
\item $\mathbb{E}''^2=0$ and thus defines a holomorphic structure
on $E$, \item $\mathbb{E}''\circ \Phi+\Phi\circ\mathbb{E}''=0$ and
so $\Phi$ is a holomorphic section, \item and $\Phi\wedge\Phi=0$.
\end{enumerate}
The statement about the dualizing module is also clear. The
duality theorem in this case is due to Simpson, \cite{Sim1}.
\endproof
\subsubsection{Generalized Higgs Algebroids} The example above is a
special case of a general construction. Let $\rho:\Al\to T_\cx X$
be a Lie algebroid and $(E,\mathbb{E})$ a module over $\Al$. Then
set $\Al_E=\Al\oplus E$ with the anchor map being the composition
$\Al\oplus E\to \Al\stackrel{\rho}{\to}T_\cx X$. Define the
bracket as
$$
[X_1+e_1,X_2+e_2]_E=[X,Y]+\mathbb{E}_{X_1}e_2-\mathbb{E}_{X_2}e_1
$$
\begin{prop}
\begin{enumerate}
\item $\Al_E$ is a Lie algebroid. \item If $\Al$ is elliptic, then
$\Al_E$ is elliptic as well. \item A module
$(H,\mathbb{H}=\mathbb{H}_0+\Phi)$ over $\Al_E$ consists of a
triple $(H,\mathbb{H}_0,\Phi)$ where $(H,\mathbb{H}_0)$ is an
$\Al$ module and $\Phi:H\to H\otimes\Al^\vee$ satisfies
$[\mathbb{H}_0,\Phi]=0$ (i.e. $\Phi$ is a morphism of $\Al$
modules) and $\Phi\wedge\Phi=0$.
\end{enumerate}
\end{prop}
\proof All of these statements follow as in the previous example.
\endproof
We call such a triple $(H,\mathbb{H}_0,\Phi)$ a Higgs bundle with
coefficients in $E$.

\subsubsection{The generalized complex Lie algebroid}
 Recall from \cite{Hit} and \cite{Gu} that an
almost generalized complex structure on a manifold $X$ is defined
by a subbundle
$$E\subseteq (TX\oplus T^\vee X)_\cx$$ satisfying
$E$ is a maximal isotropic complex subbundle $E\subset (TX\oplus
TX^\vee)_\cx$ such that $E\cap\overline{E}=\{0\}$. The isotropic
condition is with respect to the bilinear form
\[
\langle X+\xi,~ Y+\eta\rangle=\frac{1}{2}(\xi(Y)+\eta(x))
\]
The almost generalized complex structure $E$ is integrable and $E$
is called a generalized complex structure if the sections of $E$,
$\Gamma (E)$ are closed under the Courant bracket. The Courant
bracket is a skew-symmetric bracket defined on smooth sections of
$(TX\oplus TX^\vee)_\cx$, given by
\begin{equation*}
[X+\xi,Y+\eta]=[X,Y]+L_X\eta-L_Y\xi-\frac{1}{2}d(i_X\eta-i_Y\xi),
\end{equation*}
where $X+\xi,Y+\eta\in \Gamma(TX\oplus TX^\vee)_\cx$. It is shown
in \cite{Hit} and \cite {Gu} how symplectic and complex manifolds
are examples of generalized complex manifolds.

\np In the case of a generalized complex structure, the projection
map $\rho:E\to T_\cx X$ defines a Lie algebroid, the Lie algebra
structure on the sections of $E$ being the Courant bracket. Note
that the Courant bracket on the full space $(TX\oplus
TX^\vee)\otimes\cx$ does not satisfy Jacobi.

\begin{prop}
$E$ is an elliptic Lie algebroid.
\end{prop}
\proof That it is a Lie algebroid is a straight forward
calculation, as in \cite{Gu}. That it is elliptic follows just as
in the case of the Higgs Lie algebroid. \endproof

\np In this case, Gualitieri \cite{Gu} calls cohesive modules generalized holomorphic vector bundles. 
There is therefore a duality theorem in this context. In
general it can not be made any more explicit than the general
duality theorem (\ref{dualitytheorem}). On the other hand, in the
special case where the generalized complex manifold is a complex
manifold $X$, $E=T^{0,1}X\oplus T^{1,0\vee}X$, \cite{Gu}. Then we
have
\begin{prop}
Let $E$ be the algebroid coming from the generalized complex
structure defined by an honest complex structure as defined above.
Then
\begin{enumerate}
\item a module over $E$ consists of the following data: $(E,\Phi)$
where $E$ is a holomorphic vector bundle on $X$, and
$\Phi\in\mbox{Hom}(E,E\otimes T^{1,0})$ is a holomorphic section
and satisfies $\Phi\wedge\Phi=0$.

\item The dualizing module $\D_E$ is $(K^{\otimes 2},0)$ the
square of the canonical bundle with the zero Higgs field $\Phi$
\end{enumerate}
\end{prop}
\proof The proof is the same as for the Higgs algebroid.
\endproof

\

\subsection{Non-commutative tori}
\subsubsection{Real noncommutative tori}\label{realnct}
We now describe noncommutative tori. We will describe
them in terms of twisted group algebras. Let $V$ be a real vector
space, and $\Lambda\subset V$ a lattice subgroup. The we can form
the group ring $\Sc^*(\Lambda)$, the Schwartz space of complex
valued functions on $\Lambda$ which decrease faster than any
polynomial. Let $B\in \Lambda^2 V^\vee$, and form the biadditive,
antisymmetric group cocycle $\sigma:\Lambda\times\Lambda\to U(1)$
by \[ \sigma(\lambda_1,\lambda_2)=e^{2\pi i
B(\lambda_1,\lambda_2)}
\]
In our computations, we will often implicitly make use of the fact
that $\sigma$ is biadditive and anti-symmetric.  Now we can form
the twisted group algebra $\A(\Lambda;\sigma)$ consisting of the
same space of functions as $\Sc^*(\Lambda)$ but where the
multiplication is defined by
\[
[\lambda_1]\circ[\lambda_2]=\sigma(\lambda_1,\lambda_2)[\lambda_1+\lambda_2]
\]
This is a $*$-algebra where
$f^*(\lambda)=\overline{f(\lambda^{-1})}$. This is one of the
standard ways to describe the (smooth version) of the
noncommutative torus.  Given $\xi\in V^\vee$ it is easy to check
that
\begin{equation}\label{eq-derivation}
\xi(f)(\lambda)=2\pi \imath\langle \xi, \lambda\rangle f(\lambda)
\end{equation}
defines a derivation on $\A(\Lambda;\sigma)$. Note that the
derivation $\xi$ is ``real" in the sense that $\xi(f^*)=-\xi(f)$.
Finally define a (de Rham) dga $\As$ by
\[
\A^\bullet(\Lambda;\sigma)=\A(\Lambda;\sigma)\otimes
\Lambda^\bullet V_\mathbb{C}\] where $V_\mathbb{C}=V\otimes
\mathbb{C}$ and the differential $d$ is defined on functions
$\phi\in \A(\Lambda;\sigma)$ by
\[
\langle df, \xi\rangle=\xi(f)
\]
for $\xi\in V_\mathbb{C}^\vee$. In other words, for $\lambda\in
\Lambda$ one has $d\lambda=2\pi\imath\lambda \otimes D(\lambda)$
where $D(\lambda)$ denotes $\lambda$ as an element of $\Lambda^1
V$. Extend $d$ to the rest of $\A^\bullet(\Lambda;\sigma)$ by
Leibniz. Note that $d^2=0$.

\begin{rem} We just want to point out that $V$ appears as the
``cotangent" space. This is a manifestation of the fact that there
is a duality going on. That is, in the case that $\sigma=1$, we
have that the dga $\As=(\A^\bullet(\Lambda;\sigma),d,0)$ is naturally
isomorphic to $(\A^\bullet(V^{\vee}/\Lambda^\vee),d ,0)$ the de Rham
algebra of the dual torus and $T^\vee_0 (V^\vee/\Lambda^\vee)$ is
naturally isomorphic to $V$. See Proposition \ref{prop-dualDGA}
for the complex version of this.
\end{rem}

\subsubsection{Complex noncommutative tori}\label{complexnct}
We are most interested in the case where our torus has a complex
structure and in defining the analogue of the Dolbeault DGA for a
noncommutative complex torus.  So now let $V$ will be a vector
space with a complex structure $J:V\to V$, $J^2=-\mathbb{1}$. Let
$g$ be the complex dimension of $V$. Set
$V_\mathbb{C}=V\otimes_\mathbb{R}\mathbb{C}$. Then $J\otimes
1:V_\mathbb{C}\to V_\mathbb{C}$ still squares to $-\mathbb{1}$ and
so $V_\mathbb{C}$ decomposes into $\imath$ and $-\imath$
eigenspaces, $V_{1,0}\oplus V_{0,1}$. The dual $V_\mathbb{C}^\vee$
also decomposes as $V_\mathbb{C}^\vee=V^{1,0}\oplus V^{0,1}$.  Let
$D':V_\mathbb{C}\otimes \mathbb{C}\to V_{1,0}$ and $D'':V\otimes
\mathbb{C}\to V_{0,1}$ denote the corresponding projections.
Explicitly
\[D'=\frac{J\otimes 1 + 1\otimes \imath}{2\imath}\] and
\[D''=\frac{-J\otimes 1 + 1\otimes \imath}{2\imath}\]
and $D=D'+D''$ where $D$ denotes the identity. This also
established a decomposition
\[\Lambda^kV_\mathbb{C}=\otimes_{p+q=\bullet}\Lambda^{p,q}V\]
where $\Lambda^{p,q}V=\Lambda^p V_{1,0}\otimes \Lambda^qV_{0,1}$.
Complex conjugation on $V_\mathbb{C}$ defines an involution and
identifies $V$ with the $v\in V_\mathbb{C}$ such that
$\overline{v}=v$.

Now let $X=V/\Lambda$, a complex torus of dimension $g$,  and
$X^\vee=\overline{V}^\vee/\Lambda^\vee$ its dual torus. Let
$B\in\Lambda^2 V^\vee$ be a real (constant) two form on $X$. Then
$B$ will decompose in to parts \[ B=B^{2,0}+B^{1,1}+B^{0,2}\]
where $B^{p,q}\in \Lambda^{p,q}V^\vee$,
$B^{0,2}=\overline{B^{2,0}}$ and $\overline{B^{1,1}}=B^{1,1}$. Now
$B^{0,2}\in \Lambda^2 V^{0,1}\cong H^{0,2}(X)$. Then it also
represents a class  $\Pi\in \Lambda^2 V^{0,1}\cong
H^0(X^\vee;\Lambda^2 T_{1,0}X)$. Let $\sigma:\Lambda \wedge
\Lambda\to U(1)$ denote the group $2$-cocycle given by
\[\sigma(\lambda_1,\lambda_2)=e^{2\pi
\imath B(\lambda_1,\lambda_2)}.\] and form as above
$\A(\Lambda;\sigma)$ the twisted group algebra based on rapidly
decreasing functions. Define the Dolbeault dga $\As$
$\A^{0,\bullet}(\Lambda;\sigma)$ to be
\[
\A(\Lambda;\sigma)\otimes \Lambda^\bullet V_{1,0}
\]
where for $\lambda\in \A(\Lambda;\sigma)$ we define
\[
\dbar\lambda=2\pi\imath \lambda\otimes D'(\lambda)\in
\A(\Lambda;\sigma)\otimes V_{1,0}
\]
We can then extend $\dbar$ to the rest of
$\A^{0,\bullet}(\Lambda;\sigma) $ by the Leibniz rule. Let us
reiterate the remarks above. Even though we are defining the
$\dbar$ operator, we are using the $(1,0)$ component of
$V_\mathbb{C}$. This is because of duality. In the case of the
trivial cocycle $\sigma$, this definition is meant to reconstruct
the Dolbeault algebra on $X^\vee$. In this case,
$\A^{0,\bullet}(X^\vee)\cong \A^\bullet(X^\vee)\otimes
\Lambda^\bullet T^{0,1}_0X^\vee$. But
\[T^{0,1}_0X^\vee = \overline{\overline{V}^\vee}^\vee\cong
V_{1,0}.\] To check the reasonableness of this definition we have
\begin{prop}\label{prop-dualDGA} If $\sigma=1$ is the trivial
cocycle, then the dga $(\A^{0,\bullet}(\Lambda;\sigma), \dbar)$ is
isomorphic to the Dolbeault dga $(\A^{0,\bullet}(X^\vee),\dbar)$.
\end{prop}

We now show that the dga $\As=(\A^{0,\bullet}(\Lambda;\sigma), \dbar,0)$ is elliptic. Let
\[
\tau:\A(\Lambda;\sigma)\to \mathbb{C}
\]
denote the continuous $\mathbb{C}$-linear functional defined by
$\tau(\sum a_\lambda\lambda)=a_0$. This is a trace, that is
$\tau(ab)=\tau(ba)$ and in many cases it is the unique normalized
trace on $\A(\Lambda;\sigma)$. (It is unique when $\sigma$ is
``irrational" enough.) We note the following lemma whose proof is
straightforward.
\begin{lem}
For any $\xi\in V^\vee$, the derivation $\xi$ defined by
\eqref{eq-derivation} has the property
\begin{equation}\label{eq-nctstokes}
\tau(\xi(f))=0
\end{equation}
for all $f\in \A(\Lambda;\sigma)$.
\end{lem}

\noindent  Equip $V_\mathbb{C}$ with an Hermitian inner product
$\langle \cdot,\cdot\rangle:V_\mathbb{C}\times
V_\mathbb{C}\to\mathbb{C}$. Let $v_1, \ldots, v_g$ and $v^1,
\ldots, v^g$ be dual orthonormal bases of $V_\mathbb{C}$ and
$V_\mathbb{C}^\vee$ respectively. Equip $V_\mathbb{C}$ with a
Hermitian structure. Then $V_{1,0}$ and $V_{0,1}$ inherit
Hermitian structures as well. Let $v'_i$ and $v''_i$
($i=1,\cdots,g$) be orthonormal bases of $V_{1,0}$ and $V_{0,1}$
respectively. We let $D=\A(\Lambda;\sigma)\otimes \Lambda^g
V_{0,1}$ with
\[
\mathbb{D}:D\to
D\otimes_{\A(\Lambda;\sigma)}\A^{0,1}(\Lambda;\sigma)
\]
defined by
\[
\mathbb{D}(f\otimes v_{\{1,\cdots,g\}}'')=\dbar(f)\otimes
v_{\{1,\cdots,g\}}''
\]
Recall that the $V_{1,0}$ is the {\em anti-holomorphic cotangent}
space of the noncommutative complex torus and $V_{0,1}$ is the
{\em holomorphic cotangent space}.

\noindent Define $\bar{*}:\A^{0,k}(\Lambda;\sigma)\to
D\otimes_{\A(\Lambda;\sigma)}\A^{0,g-k}$ by
\[
\bar{*}(f\otimes v_I')=v_{\{1,\cdots,g\}}''\otimes f^* \otimes
v_{I^c}'
\]

\

\noindent Now note that
$D\otimes_{\A(\Lambda;\sigma)}\A^{0,g}(\Lambda;\sigma)\cong
\A^{2g}(\Lambda;\sigma)$ and so we define
\[
\int
:D\otimes_{\A(\Lambda;\sigma)}\A^{0,g}(\Lambda;\sigma)\to\mathbb{C}
\]
Define \begin{equation} \int a_\lambda\lambda\otimes
v'_I\wedge v''_I=\left\{\begin{array}{ll} 0 & \mbox{ if } I\ne \{1,\cdots, g\}\\
\tau(a_\lambda) & \mbox{ if } I=\{1,\cdots, g\}
\end{array}\right\}
\end{equation}
The following lemma is
trivial to verify.
\begin{lem}
For all $x\in
D\otimes_{\A(\Lambda;\sigma)}\A^{0,g}(\Lambda;\sigma)$ we have
\[
\int \mathbb{D}(x)=0
\]
\end{lem}
\begin{thm}(Serre duality for complex noncommutative tori)
\begin{enumerate}
\item The dga 
$\As=(\A^{0,\bullet}(\Lambda;\sigma),\dbar,0$ is elliptic with dualizing module $(D,\mathbb{D}),\bar{*},\int)$. \item On the category
$\Pc_{\A^{0,\bullet}(\Lambda;\sigma) }$ there is a Serre functor
defined by
\[
(E^\bullet,\mathbb{E})\mapsto
(E\otimes_{\A(\Lambda;\sigma)}D,\mathbb{E}\#\mathbb{D})
\]
\item In the case when $\sigma=1$, the Serre functor coincides
with the usual Serre functor on $X^\vee$ using the isomorphism
described in Proposition \ref{prop-dualDGA}
\end{enumerate}
\end{thm}

\bibliographystyle{amsalpha}

\end{document}